\documentclass[11pt,a4paper]{article}
\usepackage[T1]{fontenc}
\usepackage[utf8]{inputenc}
\usepackage{lmodern}
\usepackage[margin=25mm]{geometry}
\usepackage{amsmath,amssymb,amsfonts,amsthm,amsopn}
\usepackage{mathrsfs,mathtools,blkarray}
\usepackage{graphicx,float,enumerate,booktabs,multirow,diagbox}
\usepackage{algorithm,algorithmic}
\usepackage{subcaption,rotating}
\usepackage{xcolor,tikz}
\usetikzlibrary{patterns}
\usepackage{xurl}
\usepackage[hidelinks]{hyperref}
\usepackage{cleveref}

\numberwithin{equation}{section}
\numberwithin{algorithm}{section}
\theoremstyle{plain}
\newtheorem{theorem}{Theorem}[section]
\newtheorem{proposition}[theorem]{Proposition}
\newtheorem{definition}[theorem]{Definition}
\newtheorem{problem}[theorem]{Problem}
\newcommand{\email}[1]{\href{mailto:#1}{#1}}
\DeclareMathOperator{\diag}{diag}
\colorlet{FColor}{orange}
\colorlet{LRColor}{cyan}

\colorlet{DenseColor}{gray!60}
\def\diagsz{1.8}
\def\restdiag{0.1}

\def\bigsz{0.65}
\def\smallsz{0.2}
\def\off{0.05}

\newcommand*{\herm}{*}
\newcommand*{\pinv}{\dagger}
\newcommand*{\expectation}{\mathbb{E}}

\newcommand*{\complex}{\mathbb{C}}
\newcommand*{\unitCircle}{\mathbb{S}}

\newcommand*{\integer}{\mathbb{Z}}
\newcommand*{\Mat}[1]{\boldsymbol{#1}}

\newcommand*{\rank}{\operatorname{rank}}
\newcommand*{\range}{\operatorname{range}}
\newcommand*{\bigO}{\mathscr{O}}
\newcommand*{\e}{\mathrm{e}}
\newcommand*{\imag}{\mathrm{i}}
\newcommand*{\compl}{\mathrm{c}}
\newcommand*{\fro}{\mathrm{F}}
\newcommand*{\diff}{\mathrm{d}}
\newcommand*{\HSS}{\mathrm{HSS}}
\newcommand*{\fast}{\mathrm{fast}}

\newcommand*{\hierarchical}{\mathcal{H}}
\newcommand*{\indexset}[1]{{\mathcal{#1}}}
\newcommand*{\tree}{\mathsf{T}}
\newcommand*{\mult}{\mathrm{mult}}

\newcommand*{\level}{\operatorname{level}}
\newcommand*{\bg}{\mathrm{big}}

\newcommand*{\solve}{\mathrm{s}}
\newcommand*{\err}{\mathrm{err}}

\newcommand*{\cB}{\mathrm{cB}}
\newcommand*{\fB}{\mathrm{fB}}
\newcommand*{\cH}{\mathrm{cH}}
\newcommand*{\fH}{\mathrm{fH}}
\newcommand*{\s}{\mathrm{s}}
\newcommand*{\direct}{\mathrm{d}}
\newcommand*{\pre}{\mathrm{pre}}
\newcommand*{\iter}{\mathrm{iter}}
\newcommand*{\exact}{\mathrm{ex}}

\title{A Superfast Direct Solver for Type-III Inverse Nonuniform Discrete Fourier Transform}
\author{
Yingzhou Li\thanks{School of Mathematical Sciences, Fudan University, 220 Handan Road, Shanghai 200433, China;
Shanghai Key Laboratory for Contemporary Applied Mathematics, Fudan University, 220 Handan Road, Shanghai 200433, China
  (\email{yingzhouli@fudan.edu.cn}).}
\and Jingyu Liu\thanks{School of Mathematical Sciences, Fudan University, 220 Handan Road, Shanghai 200433, China
  (\email{jyliu22@m.fudan.edu.cn}).}
}
\date{}
\hypersetup{pdftitle={A Superfast Direct Solver for Type-III Inverse Nonuniform Discrete Fourier Transform},pdfauthor={Yingzhou Li and Jingyu Liu}}

\begin{document}

\maketitle

\begin{abstract}
Various fields in scientific computing require the use of the nonuniform discrete Fourier transform~(NUDFT) and its inversion.
In~[H.~Wilber, E.~N.~Epperly, and A.~H.~Barnett, \textit{SIAM J. Sci. Comput.}, 47 (2025), pp.~A1702--A1732], the authors proposed a direct inversion method for the type-II NUDFT, which achieves quasi-linear complexity by approximating the type-II NUDFT matrix as a product of a hierarchical semiseparable~(HSS) matrix and a DFT matrix.
In this article, we propose a novel superfast direct inversion method for type-III NUDFT.
The proposed method approximates the type-III NUDFT matrix as a product of a type-II NUDFT matrix and an HSS matrix.
This decomposition enables both the forward application and the backward inversion to be accomplished with quasi-linear complexity.
The fast inversion can serve as a high-accuracy direct solver or as an efficient preconditioner.
Additionally, we provide an error bound for the approximation when the sample points are independent and identically distributed uniform random variables and the frequencies are perturbed continuous integers.
Numerical results are presented to verify the relevant theoretical properties and demonstrate the efficiency of the proposed methods.
\end{abstract}

\par

\par\medskip

\section{Introduction} \label{sec:introduction}

In this paper, a superfast direct solver is proposed for the \textit{nonuniform discrete Fourier transform}~(NUDFT).
In one dimension, the NUDFT \textit{forward} problem has the following form:
\begin{equation} \label{eq:nudft}
    f_{j} = \sum_{k = 0}^{N - 1} c_{k} \e^{-2 \pi \imag  \omega_{k} x_{j}}, \quad  0 \leq j \leq M - 1,
\end{equation}
where \(\{x_{j}\} \subset [0, 1)\) are \textit{sample points}, \(\{\omega_{k}\} \subset [-1 / 2, N - 1 / 2)\) are \textit{frequencies}, \(\{c_{k}\} \subset \complex\) are \textit{coefficients}, and \(\{f_{j}\} \subset \complex\) are \textit{target values}.
We assume throughout the paper that \(M \geq N\).
Since the forward problem is a matrix-vector multiplication, direct computation requires \(\bigO(M N)\) operations.
If \(M = N\), the sample points are equispaced and the frequencies are continuous integers, i.e., \(x_{j} = j / M\) and \(\omega_{k} = k\), the NUDFT becomes uniform, and~\eqref{eq:nudft} can be computed using the \textit{fast Fourier transform}~(FFT) algorithm in \(\bigO(N \log N)\) operations~\cite{Cooley_Tukey_1965}.
Unfortunately, the FFT cannot be directly applied when either the sample points or the frequencies are nonuniform.
In these cases, the NUDFT problem can be categorized into three types according to the nonuniformity~\cite{Dutt_Rokhlin_1995, Ruiz-Antolín_Townsend_2018}:
\begin{itemize}
    \item Type-I: Uniform sample points and nonuniform frequencies, i.e., \(x_{j} = j / M\);
    \item Type-II: Nonuniform sample points and uniform frequencies, i.e., \(\omega_{k} = k\);
    \item Type-III: Nonuniform sample points and nonuniform frequencies.
\end{itemize}
Various algorithms have been developed to address the NUDFT forward problem with a complexity of \(\bigO(M + N + N \log N)\) while maintaining the desired accuracy.
These algorithms are commonly referred to as the \textit{nonuniform fast Fourier transform}~(NUFFT).
Some methods approximate the NUDFT matrix by the Hadamard product of the discrete Fourier transform~(DFT) matrix and a low-rank matrix~\cite{Anderson_Dahleh_1996, Ruiz-Antolín_Townsend_2018}, while others utilize convolution with a smooth, nearly locally supported function combined with the FFT on an oversampled grid~\cite{Barnett_Magland_Klinteberg_2019, Dutt_Rokhlin_1993, Dutt_Rokhlin_1995, Greengard_Lee_2004, Lee_Greengard_2005,Plonka_Potts_Steidl_Tasche_2023, Potts_Steidl_Tasche_2001}.
It is worthwhile to point out that typically in the NUFFT algorithms, the linear term \(\bigO(M + N)\) comes from the number of sample points and frequencies, while the quasi-linear term \(\bigO(N \log N)\) is from the fact that the frequencies are distributed in an interval of length \(\bigO(N)\).
For more general problems such as the Fourier integral operator~(FIO), \textit{butterfly algorithms} have been designed for efficient forward computation~\cite{Candes_Demanet_Ying_2009, Li_Yang_Martin_Ho_Ying_2015} at a complexity of \(\bigO(N \log N )\) when \(M = N\).

In this article, we are mainly interested in the \textit{inverse} NUDFT problem.
\begin{problem}[Inverse NUDFT (INUDFT)] \label{problem:inudft}
    Given sample points \(\{x_{j}\}\) and frequencies \(\{\omega_{k}\}\), determine the coefficients \(\{c_{k}\}\) from target values \(\{f_j\}\).
\end{problem}
Algebraically, this problem can be modeled as a linear least-squares problem associated with the NUDFT matrix \(\Mat{A}\), whose entries are given by \(\Mat{A}(j, k) = \e^{-2 \pi \imag x_{j} \omega_{k}}\).
Unlike the uniform case, the pseudoinverse of the NUDFT matrix is not simply equal to its adjoint~(up to a scaling factor), meaning that the fast inverse algorithm differs from that used for the forward problem.

Another alternative perspective on the inverse problem is to view it as a function approximation task.
Specifically, given values \(\{f(x_{j})\}\) of a function \(f \colon [0, 1] \to \complex\) at sample points \(\{x_{j}\}\), we aim to determine the coefficients \(\{c_{k}\}\) such that the expansion \(\sum_{k = 0}^{N - 1} c_{k} \e^{-2 \pi \imag \omega_{k} x}\) provides a good approximation of \(f\).
According to the Kadec-\(1/4\) theorem~\cite{Young_1981_book}, when frequencies are perturbed integers satisfying \(\sup_{k \in \integer} |\omega_{k} - k| \leq \alpha \) with \(0 \leq \alpha < 1 / 4\),
the nonharmonic Fourier modes \(\{\e^{-2 \pi \imag \omega_{k} x}\}_{k \in \integer}\) form a Riesz basis of \(L^{2}[0, 1]\).
Hence the INUDFT problem can be interpreted as finding the coefficients of a given function on some nonharmonic Fourier basis.
For more applications of the type-III NUDFT such as the heat flow and the magnetic resonance imaging~(MRI), we refer the reader to~\cite{Lee_Greengard_2005}.

Typically, the INUDFT problem is solved using iterative methods~\cite{Dutt_Rokhlin_1993, Feichtinger_Grochenig_Strohmer_1995, Ruiz-Antolín_Townsend_2018} due to the efficiency of its forward computation.
Direct methods have also been developed, as noted in~\cite{Kircheis_Potts_2019,Kircheis_Potts_2023}.
The authors design a direct solver for the type-II problem by finding an approximation of the form \(\Mat{A} \Mat{D}^{\herm} \Mat{F}^{\herm} \Mat{S} \approx N \Mat{I}_{M}\) where the sparse matrix \(\Mat{S}\) is the only matrix to be determined.
However, these direct methods are for the type-I or type-II problems and cannot be easily extended to the type-III problem.
Recently, in~\cite{Wilber_Epperly_Barnett_2025}, the authors proposed a direct solver for the type-II INUDFT problem.
A key observation they use is that type-II NUDFT matrix \(\Mat{A}\) is a Vandermonde matrix and thus satisfies a Sylvester matrix equation with a low-rank right-hand side.
Consequently, the transformed matrix \(\Mat{\tilde{A}} = \Mat{A} \Mat{F}^{-1}\) can be compressed into a \textit{hierarchically semiseparable}~(HSS) matrix~\cite{Xia_Chandrasekaran_Gu_Li_2010}, where \(\Mat{F} \in \complex^{N \times N}\) is the DFT matrix defined by \(\Mat{F}(j, k) = \e^{-2 \pi \imag j k / N}\).
Their direct inversion method is achieved by applying the HSS least-squares solver~\cite{Xi_Xia_Cauley_Balakrishnan_2014} and the \textit{inverse fast Fourier transform}~(iFFT).
The proposed method could be naturally extended to the type-I INUDFT problem but not to the type-III INUDFT problem, since the type-III NUDFT matrix is no longer a Vandermonde matrix and no such Sylvester matrix equation can be found.

\subsection{Contributions} \label{subsec:contributions}

The main purpose of this paper is to develop a direct solver for the type-III INUDFT problem.
We only consider the case where the frequencies are the perturbed continuous integers.
By projecting the frequencies onto the integers, the type-III NUDFT matrix is approximated as \(\Mat{A} \approx \Mat{B} \Mat{H}\), where \(\Mat{B} \in \complex^{M \times N}\) is a type-II NUDFT matrix and \(\Mat{H} \in \complex^{N \times N}\) is a dense matrix that can be compressed into an HSS matrix.
We provide a theoretical error bound for this approximation under the assumption that the sample points are independent and identically distributed~(i.i.d.)~uniform random variables.

By applying relevant fast inversion algorithms to \(\Mat{B}\) and \(\Mat{H}\), the least-squares problem associated with \(\Mat{A}\) can be solved rapidly.
To be more specific, the usage of this superfast direct solver can be divided into four stages: preliminary, construction, factorization, and solution.
In the preliminary stage, we build a fast structure \(\Mat{B}_{\fast}\) for the type-II NUDFT matrix \(\Mat{B}\).
In the construction stage, the matrix \(\Mat{B}_{\fast}^{\pinv} \Mat{A}\) is compressed into an HSS matrix, denoted as \(\Mat{H}_{\HSS}\).
This compression is performed efficiently using an algorithm based on random sampling.
In the factorization stage, the HSS matrix \(\Mat{H}_{\HSS}\) undergoes the URV factorization~\cite{Xi_Xia_Cauley_Balakrishnan_2014}.
Once the factorization stage finishes, we use \(\Mat{A}_{\fast}^{\pinv} =  \Mat{H}_{\HSS}^{-1} \Mat{B}_{\fast}^{\pinv}\) as a direct solver for the least-squares problem.
Additionally, when considering the normal equation of~\eqref{eq:nudft}, our solver can also serve as an efficient preconditioner in preconditioned conjugate gradient~(PCG) methods.
Numerical experiments are conducted to investigate the properties of the approximation and to demonstrate the efficiency of the proposed method.

\subsection{Organization} \label{subsec:organization}

The rest of this paper is organized as follows.
In section~\ref{sec:preliminaries}, we review the related definitions and algorithms related to HSS matrices, as well as the type-II INUDFT solver~\cite{Wilber_Epperly_Barnett_2025}, which serves as a preliminary for the subsequent sections.
Section~\ref{sec:inudft3} introduces our novel method on the type-III NUDFT problem.
We discuss first its derivation and related theoretical results, followed by the practical implementation of the algorithm.
Numerical results are presented in section~\ref{sec:numerical_results} to illustrate the efficiency of our methods.
Finally, section~\ref{sec:conclusions_future_work} concludes the paper with a discussion on future directions.

\section{Preliminaries} \label{sec:preliminaries}

\subsection{HSS matrices} \label{subsec:hss}

HSS matrices are a special class of \(\hierarchical^2\)-matrices~\cite{Borm_Grasedyck_Hackbusch_2003, Hackbusch_1999, Hackbusch_Khoromskij_Sauter_2000} which have linear complexity for storage and numerical linear algebra operations, such as matrix-vector products and solving linear systems~\cite{Xia_Chandrasekaran_Gu_Li_2010}.
The definition of an HSS matrix relies on a tree structure~(HSS tree), which is also included in the following definition.
\begin{definition}[HSS tree and HSS matrix~\cite{Xi_Xia_Cauley_Balakrishnan_2014}] \label{def:hss}
    An \(M \times N\) matrix \(\Mat{A}\) is called an HSS matrix with the associated HSS tree \(\tree\) if the following conditions hold:
    \begin{itemize}
        \item \(\tree\) is a full binary tree with root node \(\rho\).
        \item There are two index sets \(\indexset{J}_{\tau}\) and \(\indexset{K}_{\tau}\) associated with each node \({\tau}\) of \(\tree\), satisfying the following conditions:
        \begin{enumerate}[(1)]
            \item \(\indexset{J}_{\rho} = \{0, \dotsc, M - 1\}\), \(\indexset{K}_{\rho} = \{0, \dotsc, N - 1\}\).
            \item For a nonleaf node \(\tau\) with children \(\alpha_{1}\) and \(\alpha_{2}\), we have \(\indexset{J}_{\tau} = \indexset{J}_{\alpha_{1}} \sqcup \indexset{J}_{\alpha_{2}}\) and \(\indexset{K}_{\tau} = \indexset{K}_{\alpha_{1}} \sqcup \indexset{K}_{\alpha_{2}}\), where the notation \(\sqcup\) means the disjoint union.
        \end{enumerate}
        \item There are matrices \(\Mat{D}_{\tau}\) and \(\Mat{U}_{\tau}\), \(\Mat{V}_{\tau}\), \(\Mat{B}_{\tau, \sigma}\) called HSS generators associated with each node \(\tau\) and its sibling \(\sigma\), satisfying the following recursion:
        \begin{enumerate}[(1)]
            \item
            If \(\tau\) is a leaf node, \(\Mat{D}_{\tau} = \Mat{A}({\indexset{J}_{\tau}, \indexset{K}_{\tau}})\) is a dense matrix.
            If \(\tau\) is a nonleaf node with children \(\alpha_{1}\) and \(\alpha_{2}\),
            \begin{equation*}
                \Mat{D}_{\tau}
                = \Mat{A}({\indexset{J}_{\tau}, \indexset{K}_{\tau}})
                =
                \begin{bmatrix}
                    \Mat{D}_{\alpha_{1}} & \Mat{U}_{\alpha_{1}}^{\bg} \Mat{B}_{\alpha_{1}, \alpha_{2}} \Mat{V}_{\alpha_{2}}^{\bg, \herm} \\
                    \Mat{U}_{\alpha_{2}}^{\bg} \Mat{B}_{\alpha_{2}, \alpha_{1}} \Mat{V}_{\alpha_{1}}^{\bg, \herm} & \Mat{D}_{\alpha_{2}} \\
                \end{bmatrix},
            \end{equation*}
            where \(\Mat{V}_{\alpha_{j}}^{\bg, \herm}\) means \((\Mat{V}_{\alpha_{j}}^{\bg})^{\herm}\).
            This notation will be used throughout the paper.
            \item If \(\tau\) is a leaf node, we define \(\Mat{U}_{\tau}^{\bg} = \Mat{U}_{\tau}\) and \(\Mat{V}_{\tau}^{\bg} = \Mat{V}_{\tau}\).
            When \(\tau\) is a nonleaf and nonroot node, the matrices \(\Mat{U}_{\tau}\) and \(\Mat{V}_{\tau}\) satisfy the following recursions:
            \begin{equation*}
                \Mat{U}_{\tau}^{\bg} =
                \begin{bmatrix}
                    \Mat{U}_{\alpha_{1}}^{\bg} & \\
                    & \Mat{U}_{\alpha_{2}}^{\bg} \\
                \end{bmatrix}
                \Mat{U}_{\tau}
                \quad \text{and} \quad
                \Mat{V}_{\tau}^{\bg} =
                \begin{bmatrix}
                    \Mat{V}_{\alpha_{1}}^{\bg} & \\
                    & \Mat{V}_{\alpha_{2}}^{\bg} \\
                \end{bmatrix}
                \Mat{V}_{\tau}.
            \end{equation*}
        \end{enumerate}
    \end{itemize}
    Given an HSS tree, we define the level of each node as follows:
    For the root node \(\rho\), we set \(\level(\rho) = 0\).
    If \(\tau\) is a nonleaf node and \(\alpha\) is a child of \(\tau\), we set \(\level(\alpha) = \level(\tau) + 1\).
    The maximum level of the tree is denoted by \(L\).
\end{definition}

A key observation regarding HSS matrices is that for any nonroot node \(\tau\), the matrices \(\Mat{U}_{\tau}^{\bg}\) and \(\Mat{V}_{\tau}^{\bg, \herm}\) are bases for the column or row spaces of the corresponding blocks~(called the \textit{HSS blocks}) \(\Mat{A}(\indexset{J}_{\tau}, \indexset{K}_{\tau}^{\compl})\) and \(\Mat{A}(\indexset{J}_{\tau}^{\compl}, \indexset{K}_{\tau})\), respectively.
Here \(\indexset{J}_{\tau}^{\compl} = \indexset{J}_{\rho} \backslash \indexset{J}_{\tau}\) and \(\indexset{K}_{\tau}^{\compl} = \indexset{K}_{\rho} \backslash \indexset{K}_{\tau}\).
This is often called the \textit{shared basis property}.
Matrices \(\Mat{U}_{\tau}^{\bg}\) and \(\Mat{V}_{\tau}^{\bg}\) are called \textit{basis matrices}.
The \textit{HSS rank} of \(\Mat{A}\) is defined as the maximum numerical rank among all the HSS blocks of \(\Mat{A}\).
It is well known that the complexity for both storage and the matrix-vector product of an HSS matrix with HSS rank \(k\) is \(\bigO\bigl(k (M + N)\bigr)\)~\cite{Chandrasekaran_Gu_Pals_2006}.
When the context is clear, we may omit the index sets \(\indexset{J}_{\tau}\) and \(\indexset{K}_{\tau}\) in the subscript, using \(\tau\) directly to represent the corresponding block.
For instance, \(\Mat{A}_{\tau, \sigma}\) and \(\Mat{A}_{\tau, \sigma^{\compl}}\) represent \(\Mat{A}(\indexset{J}_{\tau}, \indexset{K}_{\sigma})\) and \(\Mat{A}(\indexset{J}_{\tau}, \indexset{K}_{\sigma}^{\compl})\), respectively.

The hierarchical structure of HSS matrices is often represented in a different form known as the \textit{telescoping factorization}~\cite{Martinsson_2011, Martinsson_2019_book}.
Specifically, for each level \(1 \leq \ell \leq L\), we define the following block diagonal matrices:
\begin{equation} \label{eq:telescoping_factorization_UV}
    \Mat{U}^{(\ell)} = \diag\bigl(\Mat{U}_{\tau} \colon \level(\tau) = \ell\bigr), \quad \Mat{V}^{(\ell)} = \diag\bigl(\Mat{V}_{\tau} \colon \level(\tau) = \ell\bigr).
\end{equation}
For each level \(0 \leq \ell \leq L - 1\), we define
\begin{equation} \label{eq:telescoping_factorization_B}
    \Mat{B}^{(\ell)} = \diag\bigl(\Mat{B}_{\tau} \colon \level(\tau) = \ell\bigr) \quad \text{with} \quad \Mat{B}_{\tau} =
    \begin{bmatrix}
        \Mat{0} & \Mat{B}_{\alpha_{1}, \alpha_{2}} \\
        \Mat{B}_{\alpha_{2}, \alpha_{1}} & \Mat{0} \\
    \end{bmatrix}.
\end{equation}
Then the HSS matrix \(\Mat{A}\) can be recursively factorized as
\begin{equation} \label{eq:hss_telescoping_factorization}
    \begin{aligned}
        \Mat{A}^{(L)} & = \Mat{U}^{(L)} \Mat{A}^{(L - 1)} \Mat{V}^{(L), \herm} + \Mat{D}^{(L)}, \\
        \Mat{A}^{(\ell)} & = \Mat{U}^{(\ell)} \Mat{A}^{(\ell - 1)} \Mat{V}^{(\ell), \herm} + \Mat{B}^{(\ell)}, \quad 1 \leq \ell \leq L - 1, \\
        \Mat{A}^{(0)} & = \Mat{B}^{(0)},
    \end{aligned}
\end{equation}
where \(\Mat{A}^{(L)} = \Mat{A}\) and \(\Mat{D}^{(L)} = \diag\bigl(\Mat{D}_{\tau} \colon \level(\tau) = L\bigr)\).

\subsection{A black-box construction algorithm for HSS matrices} \label{subsec:hbs_black_box_construction}

In this section we review the randomized algorithm~\cite{Levitt_Martinsson_2024} for computing the HSS compression of a given matrix \(\Mat{A}\), under the limitation that interactions with the matrix can only be accessed through applying it and its conjugate transpose on vectors.
This scenario is often referred to as the \textit{black-box} setting~\cite{Lin_Lu_Ying_2011}.
The HSS compression of a matrix \(\Mat{A}\) is constructed from the information contained in \(\{\Mat{\Omega}, \Mat{Y} = \Mat{A} \Mat{\Omega}, \Mat{\Psi}, \Mat{Z} = \Mat{A}^{\herm} \Mat{\Psi}\}\), where \(\Mat{\Omega}\) and \(\Mat{\Psi}\) are two random \textit{test matrices} and \(\Mat{Y}\) and \(\Mat{Z}\) are corresponding \textit{sample matrices}.
We emphasize that compared to the algorithm in~\cite{Levitt_Martinsson_2024}, our algorithm needs a postprocessing step to enforce the structure of \(\Mat{B}_{\tau}\) as in~\eqref{eq:telescoping_factorization_B}, which makes it slightly different from the original one.

Suppose \(\Mat{A} \in \complex^{N \times N}\) is an HSS matrix with HSS rank \(k\) and leaf node size \(m\).
Let \(r\) be a number slightly larger than \(k\), e.g., \(r = k + 5\), and let \(\Mat{\Omega}\) and \(\Mat{\Psi}\) be \(N \times s\) Gaussian test matrices, where the number of samples \(s \geq \max(r + m, 3 r)\).
The construction can be divided into two steps:
First, we conduct a bottom-up process to compute all basis matrices \(\Mat{U}\) and \(\Mat{V}\).
Some auxiliary matrices are also computed in this step.
Next, a top-down process is employed to complete the construction by assigning matrices \(\Mat{B}\) and \(\Mat{D}\).

We begin by describing how to compute matrices \(\Mat{U}_{\tau}\) and \(\Mat{V}_{\tau}\) for a leaf node \(\tau\) in level \(\ell = L\).
By considering the row block corresponding to \(\tau\) in the equation \(\Mat{Y} = \Mat{A} \Mat{\Omega}\), we have
\begin{equation} \label{eq:hss_random_sampling}
    \Mat{Y}_{\tau} = \Mat{D}_{\tau} \Mat{\Omega}_{\tau} + \Mat{A}_{\tau, \tau^{\compl}} \Mat{\Omega}_{\tau^{\compl}}.
\end{equation}
The key observation is that if there exists an orthonormal matrix \(\Mat{P}_{\tau} \in \complex^{s \times r}\) that belongs to the null space of \(\Mat{\Omega}_{\tau}\), such that \(\Mat{\Omega}_{\tau} \Mat{P}_{\tau} = \Mat{0}\), multiplying both sides of~\eqref{eq:hss_random_sampling} by \(\Mat{P}_{\tau}\) on the right yields
\begin{equation}\label{eq:hss_random_sampling_P}
    \Mat{Y}_{\tau} \Mat{P}_{\tau} = \Mat{A}_{\tau, \tau^{\compl}} \Mat{\Omega}_{\tau^{\compl}} \Mat{P}_{\tau}.
\end{equation}
Such a matrix \(\Mat{P}_{\tau}\) can be obtained by first performing a QR factorization with column pivoting~(QRCP) on \(\Mat{\Omega}_{\tau}^{\herm}\), then taking the last \(r\) columns from the Q-factor.
Note that in~\eqref{eq:hss_random_sampling_P} the contribution from \(\Mat{A}_{\tau, \tau}\) is removed.
Then \(\Mat{Y}_{\tau} \Mat{P}_{\tau}\) can be viewed as a sampling of \(\Mat{A}_{\tau, \tau^{\compl}}\) with test matrix \(\Mat{\Omega}_{\tau^{\compl}} \Mat{P}_{\tau}\).
Consequently, the basis matrix \(\Mat{U}_{\tau}\) can be constructed by taking the first \(r\) columns from the Q-factor in the QRCP of \(\Mat{Y}_{\tau} \Mat{P}_{\tau}\).
A similar procedure using \(\Mat{\Psi}\) and \(\Mat{Z}\) yields the matrix \(\Mat{V}_{\tau}\).

With matrices \(\Mat{U}_{\tau}\) and \(\Mat{V}_{\tau}\) of \(\tau\) in hand, the diagonal block can be decomposed into two parts:
\begin{equation} \label{eq:hss_D_decompose}
    \Mat{D}_{\tau}
    = \Mat{U}_{\tau} (\Mat{U}_{\tau}^{\herm} \Mat{D}_{\tau} \Mat{V}_{\tau}) \Mat{V}_{\tau}^{\herm}
    + (\Mat{D}_{\tau} - \Mat{U}_{\tau} \Mat{U}_{\tau}^{\herm} \Mat{D}_{\tau} \Mat{V}_{\tau} \Mat{V}_{\tau}^{\herm})
    \eqcolon \Mat{U}_{\tau} \Mat{\hat{D}}_{\tau} \Mat{V}_{\tau}^{\herm} + \Mat{\check{D}}_{\tau}.
\end{equation}
The first term represents the projection of \(\Mat{D}_{\tau}\) onto basis matrices, while the second term is the remaining part.
We will address the computation of \(\Mat{\hat{D}}_{\tau}\) later in the top-down process and focus on construction of \(\Mat{\check{D}}_{\tau}\) first.
By rewriting~\eqref{eq:hss_D_decompose} as
\begin{equation*}
    \Mat{\check{D}}_{\tau}
    = \Mat{D}_{\tau}
    - \Mat{U}_{\tau} \Mat{U}_{\tau}^{\herm} \Mat{D}_{\tau} \Mat{V}_{\tau} \Mat{V}_{\tau}^{\herm}
    = (\Mat{I} - \Mat{U}_{\tau} \Mat{U}_{\tau}^{\herm}) \Mat{D}_{\tau} + \Mat{U}_{\tau} \Mat{U}_{\tau}^{\herm} \Mat{D}_{\tau} (\Mat{I} - \Mat{V}_{\tau} \Mat{V}_{\tau}^{\herm}),
\end{equation*}
it suffices to find a way to obtain \((\Mat{I} - \Mat{U}_{\tau} \Mat{U}_{\tau}^{\herm}) \Mat{D}_{\tau}\) and \(\Mat{D}_{\tau} (\Mat{I} - \Mat{V}_{\tau} \Mat{V}_{\tau}^{\herm})\).
Multiplying \((\Mat{I} - \Mat{U}_{\tau} \Mat{U}_{\tau}^{\herm})\) on the left of~\eqref{eq:hss_random_sampling} and using the fact that \(\Mat{U}_{\tau}\) forms the column basis of \(\Mat{A}_{\tau, \tau^{\compl}}\), we obtain
\begin{equation*}
    (\Mat{I} - \Mat{U}_{\tau} \Mat{U}_{\tau}^{\herm}) \Mat{Y}_{\tau} = (\Mat{I} - \Mat{U}_{\tau} \Mat{U}_{\tau}^{\herm}) \Mat{D}_{\tau} \Mat{\Omega}_{\tau}.
\end{equation*}
Therefore, \((\Mat{I} - \Mat{U}_{\tau} \Mat{U}_{\tau}^{\herm}) \Mat{D}_{\tau} =(\Mat{I} - \Mat{U}_{\tau} \Mat{U}_{\tau}^{\herm}) \Mat{Y}_{\tau} \Mat{\Omega}_{\tau}^{\pinv}\) can be recovered through a least-squares solving.
A similar process yields
\(\Mat{D}_{\tau} (\Mat{I} - \Mat{V}_{\tau} \Mat{V}_{\tau}^{\herm})
= \Bigl((\Mat{I} - \Mat{V}_{\tau} \Mat{V}_{\tau}^{\herm}) \Mat{Z}_{\tau} \Mat{\Psi}_{\tau}^{\pinv}\Bigr)^{\herm}\).
Combining these results, \(\Mat{\check{D}}_{\tau}\) is obtained through
\begin{equation*}
    \Mat{\check{D}}_{\tau}
    = (\Mat{I} - \Mat{U}_{\tau} \Mat{U}_{\tau}^{\herm}) \Mat{Y}_{\tau} \Mat{\Omega}_{\tau}^{\pinv}
    + \Mat{U}_{\tau} \Mat{U}_{\tau}^{\herm} \Bigl((\Mat{I} - \Mat{V}_{\tau} \Mat{V}_{\tau}^{\herm}) \Mat{Z}_{\tau} \Mat{\Psi}_{\tau}^{\pinv}\Bigr)^{\herm}.
\end{equation*}

Suppose \(0 \leq \ell \leq L - 1\) and we have computed the basis matrices for all nodes in level \(\ell + 1\).
By extracting these basis matrices, we obtain the following recursive decomposition of \(\Mat{A}^{(\ell + 1)}\) starting from \(\Mat{A}^{(L)} = \Mat{A}\):
\begin{equation} \label{eq:hss_decompose_ell}
    \Mat{A}^{(\ell + 1)} = \Mat{U}^{(\ell + 1)} \Mat{A}^{(\ell)} \Mat{V}^{(\ell + 1), \herm} + \Mat{\check{D}}^{(\ell + 1)},
\end{equation}
where
\begin{equation*}
    \begin{aligned}
        \Mat{U}^{(\ell + 1)} & = \diag\bigl(\Mat{U}_{\tau} \colon \level(\tau) = \ell + 1\bigr), \\
        \Mat{V}^{(\ell + 1)} & = \diag\bigl(\Mat{V}_{\tau} \colon \level(\tau) = \ell + 1\bigr), \\
        \Mat{\check{D}}^{(\ell + 1)} & = \diag\bigl(\Mat{\check{D}}_{\tau}\colon \level(\tau) = \ell + 1\bigr).
    \end{aligned}
\end{equation*}
Therefore, \(\Mat{A}^{(\ell)} = \Mat{U}^{(\ell + 1), \herm} \Mat{A}^{(\ell + 1)} \Mat{V}^{(\ell + 1)}\).
If we ``eliminate'' all nodes in level \(\ell + 1\), it can be verified that \(\Mat{A}^{(\ell)}\) remains an HSS matrix with maximum level \(\ell\).
The only requirement is to update the generators at level \(\ell\).
Specifically, for every node \(\tau\) in level \(\ell\) with children \(\alpha_{1}\) and \(\alpha_{2}\), \(\tau\) becomes a leaf node of the new HSS tree, with corresponding generators updated by
\begin{equation} \label{eq:hss_merge}
    \Mat{D}_{\tau}
    \gets
    \begin{bmatrix}
        \Mat{\hat{D}}_{\alpha_{1}} & \Mat{B}_{\alpha_{1}, \alpha_{2}} \\
        \Mat{B}_{\alpha_{2}, \alpha_{1}} & \Mat{\hat{D}}_{\alpha_{2}} \\
    \end{bmatrix}.
\end{equation}

Let \(\Mat{\Omega}^{(L)} = \Mat{\Omega}\) and \(\Mat{Y}^{(L)} = \Mat{Y}\).
By multiplying~\eqref{eq:hss_decompose_ell} with \(\Mat{\Omega}^{(\ell + 1)}\) and rearranging the equation, we obtain
\begin{equation*}
    \Mat{U}^{(\ell + 1), \herm}\bigl(\Mat{Y}^{(\ell + 1)} - \Mat{\check{D}}^{(\ell + 1)} \Mat{\Omega}^{(\ell + 1)}\bigr) = \Mat{A}^{(\ell)} \bigl(\Mat{V}^{(\ell + 1), \herm} \Mat{\Omega}^{(\ell + 1)}\bigr).
\end{equation*}
This can be interpreted as a sampling of \(\Mat{A}^{(\ell)}\).
More specifically, let \(\tau\) be a node and \(\alpha_{1}\) and \(\alpha_{2}\) be its children.
The matrices \(\Mat{\Omega}_{\tau}^{(\ell)}\) and \(\Mat{Y}_{\tau}^{(\ell)}\) are defined as
\begin{equation*}
    \Mat{\Omega}_{\tau}^{(\ell)} =
    \begin{bmatrix}
        \Mat{V}_{\alpha_{1}}^{\herm} \Mat{\Omega}_{\alpha_{1}}^{(\ell + 1)} \\
        \Mat{V}_{\alpha_{2}}^{\herm} \Mat{\Omega}_{\alpha_{2}}^{(\ell + 1)} \\
    \end{bmatrix} \quad
    \text{and} \quad
    \Mat{Y}_{\tau}^{(\ell)} =
    \begin{bmatrix}
        \Mat{U}_{\alpha_{1}}^{\herm} \bigl(\Mat{Y}_{\alpha_{1}}^{(\ell + 1)} - \Mat{\check{D}}_{\alpha_{1}} \Mat{\Omega}_{\alpha_{1}}^{(\ell + 1)}\bigr) \\
        \Mat{U}_{\alpha_{2}}^{\herm} \bigl(\Mat{Y}_{\alpha_{2}}^{(\ell + 1)} - \Mat{\check{D}}_{\alpha_{2}} \Mat{\Omega}_{\alpha_{2}}^{(\ell + 1)}\bigr) \\
    \end{bmatrix}.
\end{equation*}
Matrices \(\Mat{\Psi}_{\tau}^{(\ell)}\) and \(\Mat{Z}_{\tau}^{(\ell)}\) are defined analogously.
By the previous discussions, we can compute \(\Mat{\check{D}}_{\tau}\), \(\Mat{U}_{\tau}\), and \(\Mat{V}_{\tau}\).
This procedure is iteratively repeated until reaching the root node.
If \(\tau\) is the root node, we compute \(\Mat{D}_{\tau} = \Mat{Y}_{\tau} \Mat{\Omega}_{\tau}^{\pinv}\) directly.

Next, a top-down procedure is performed to construct matrices \(\Mat{B}\) and \(\Mat{D}\).
By comparing the factorization~\eqref{eq:hss_decompose_ell} with~\eqref{eq:hss_telescoping_factorization}, the only difference is that the matrix \(\Mat{B}_{\tau}\) in~\eqref{eq:telescoping_factorization_B} has a structure in which the diagonal blocks corresponding to its children are zero.
Therefore, the key idea here is to ``push'' the diagonal part into its children blocks until we reach the leaf nodes.

For a level \(0 \leq \ell \leq L - 1\), suppose \(\tau\) is a nonleaf node on level \(\ell\) with children \(\alpha_{1}\) and \(\alpha_{2}\).
We define the matrices \(\Mat{D}_{\tau}^{(\ell)}\), \(\Mat{U}_{\tau}^{(\ell + 1)}\), and \(\Mat{V}_{\tau}^{(\ell + 1)}\) to be the submatrices of \(\Mat{D}^{(\ell)}\), \(\Mat{U}^{(\ell + 1)}\), and \(\Mat{V}^{(\ell + 1)}\) corresponding to the node \(\tau\), respectively, and let
\begin{equation*}
    \Mat{D}_{\tau}^{(\ell)}
    =
    \begin{bmatrix}
        \Mat{D}_{\tau; \alpha_{1}, \alpha_{1}} & \Mat{D}_{\tau; \alpha_{1}, \alpha_{2}} \\
        \Mat{D}_{\tau; \alpha_{2}, \alpha_{1}} & \Mat{D}_{\tau; \alpha_{2}, \alpha_{2}} \\
    \end{bmatrix}.
\end{equation*}
The corresponding block with respect to \(\tau\) in the factorization~\eqref{eq:hss_decompose_ell} can be expressed as
\begin{equation*}
    \begin{aligned}
        \Mat{A}_{\tau}^{(\ell + 1)}
        & = \Mat{U}_{\tau}^{(\ell + 1)} \Mat{A}_{\tau}^{(\ell)} \Mat{V}_{\tau}^{(\ell + 1), \herm} + \Mat{\check{D}}_{\tau}^{(\ell + 1)} \\
        & = \Mat{U}_{\tau}^{(\ell + 1)}
        \biggl(\Mat{U}_{\tau}^{(\ell)}
        \Mat{A}_{\tau}^{(\ell - 1)}
        \Mat{V}_{\tau}^{(\ell), \herm}
        +
        \Mat{D}_{\tau}^{(\ell)}\biggr)
        \Mat{V}_{\tau}^{(\ell + 1) \herm}
        +
        \Mat{\check{D}}_{\tau}^{(\ell + 1)} \\
        & =
        \Mat{U}_{\tau}^{(\ell + 1)} \Mat{U}_{\tau}^{(\ell)}
        \Mat{A}_{\tau}^{(\ell - 1)}
        \Mat{V}_{\tau}^{(\ell), \herm} \Mat{V}_{\tau}^{(\ell + 1) \herm}
        +
        \Mat{U}_{\tau}^{(\ell + 1)}
        \Mat{D}_{\tau}^{(\ell)}
        \Mat{V}_{\tau}^{(\ell + 1) \herm}
        +
        \Mat{\check{D}}_{\tau}^{(\ell + 1)}, \\
    \end{aligned}
\end{equation*}
where \(\Mat{A}_{\tau}^{(\ell - 1)}\) has been updated.
The second and third parts in the last equation can be merged together as
\begin{equation*}
    \begin{aligned}
        &
        \begin{bmatrix}
            \Mat{U}_{\alpha_{1}} & \\
            & \Mat{U}_{\alpha_{2}} \\
        \end{bmatrix}
        \begin{bmatrix}
            \Mat{D}_{\tau; \alpha_{1}, \alpha_{1}} & \Mat{D}_{\tau; \alpha_{1}, \alpha_{2}} \\
            \Mat{D}_{\tau; \alpha_{2}, \alpha_{1}} & \Mat{D}_{\tau; \alpha_{2}, \alpha_{2}} \\
        \end{bmatrix}
        \begin{bmatrix}
            \Mat{V}_{\alpha_{1}} &  \\
            & \Mat{V}_{\alpha_{2}} \\
        \end{bmatrix}^{\herm}
        +
        \begin{bmatrix}
            \Mat{\check{D}}_{\alpha_{1}} & \\
            & \Mat{\check{D}}_{\alpha_{2}} \\
        \end{bmatrix} \\
        & \quad =
        \begin{bmatrix}
            \Mat{U}_{\alpha_{1}} & \\
            & \Mat{U}_{\alpha_{2}} \\
        \end{bmatrix}
        \begin{bmatrix}
            & \Mat{B}_{\alpha_{1}, \alpha_{2}} \\
            \Mat{B}_{\alpha_{2}, \alpha_{1}} & \\
        \end{bmatrix}
        \begin{bmatrix}
            \Mat{V}_{\alpha_{1}} &  \\
            & \Mat{V}_{\alpha_{2}} \\
        \end{bmatrix}^{\herm}
        +
        \begin{bmatrix}
            \Mat{D}_{\alpha_{1}} & \\
            & \Mat{D}_{\alpha_{2}} \\
        \end{bmatrix},
    \end{aligned}
\end{equation*}
where \(\Mat{B}_{\alpha_{i}, \alpha_{j}} = \Mat{D}_{\tau; \alpha_{i}, \alpha_{j}}\) and  \(\Mat{D}_{\alpha} = \Mat{\check{D}}_{\alpha} + \Mat{U}_{\alpha} \Mat{D}_{\tau; \alpha, \alpha} \Mat{V}_{\alpha}^{\herm}\).
Note that the matrix \(\Mat{D}_{\tau; \alpha, \alpha} \) is actually the matrix \(\Mat{\hat{D}}_{\alpha}\) defined in~\eqref{eq:hss_D_decompose}.
The process is repeated until reaching the leaf nodes.
When \(\tau\) is a leaf node, the wanted matrix \(\Mat{D}_{\tau}\) is obtained after the update, and the process is finished.

The algorithm requires \(\bigO(k)\) sample points and has a complexity of \(\bigO\bigl(k^{2} (M + N) + k T_{\mult}\bigr)\), where \(T_{\mult}\) is the complexity of applying \(\Mat{A}\) and \(\Mat{A}^{\herm}\) to a vector.
Therefore, it is particularly suitable for problems where fast computations of matrix-vector products are available.

\subsection{A superfast least-squares solver for HSS matrices} \label{subsec:ls_solver_hss}

In this section, we review the HSS least-squares solver proposed in \cite{Wilber_Epperly_Barnett_2025,Xi_Xia_Cauley_Balakrishnan_2014}.
The problem we consider is defined as \(\min_{\Mat{c} \in \complex^{N}} \|\Mat{A} \Mat{c} - \Mat{f}\|\) where \(\Mat{A} \in \complex^{M \times N}\) is an HSS matrix associated with an HSS tree whose maximum level is \(L\).
The key algorithm is the URV factorization, which can be viewed as a generalized QR factorization or the transpose version of the ULV factorization~\cite{Chandrasekaran_Gu_Pals_2006, Xia_Chandrasekaran_Gu_Li_2010}.
Generally speaking, it factorizes an HSS matrix into a product of block diagonal unitary matrices and block upper-triangular matrices, where each component in the product is easy to invert.
The computation of the factorization usually costs \(\bigO\bigl(k^{2} (M + N)\bigr)\) operations, where \(k\) is the HSS rank, and once it is finished,  \(\Mat{A}^{\pinv}\) can be efficiently applied on a vector with \(\bigO\bigl(k (M + N)\bigr)\) operations.

We give a brief description of the algorithm.
For convenience we assume that all HSS blocks of \(\Mat{A}\) have the exact rank \(k\).
The factorization begins from each leaf node.
Suppose the numbers of rows and columns of a leaf node \(\tau\) are \(m_{\tau}\) and \(n_{\tau}\), respectively.
If \(m_{\tau} \gg n_{\tau}\), a \textit{size reduction} step is taken to decrease the row size of the problem.
This introduces zeros into \(\Mat{D}_{\tau}\) and \(\Mat{U}_{\tau}\) through a QR decomposition
\begin{equation*}
    \begin{bmatrix}
        \Mat{U}_{\tau} & \Mat{D}_{\tau}
    \end{bmatrix}
    = \Mat{\Omega}_{\tau}
    \begin{bmatrix}
        \Mat{R}_{1, 1} & \Mat{R}_{1, 2} \\
        \Mat{0} & \Mat{R}_{2, 2} \\
        \hline
        \Mat{0} & \Mat{0}
    \end{bmatrix}
    = \Mat{\Omega}_{\tau}
    \begin{bmatrix}
        \Mat{\check{U}}_{\tau} & \Mat{\check{D}}_{\tau} \\
        \Mat{0} & \Mat{0} \\
    \end{bmatrix},
\end{equation*}
where \(\Mat{\Omega}_{\tau} \in \complex^{m_{\tau} \times m_{\tau}}\) is unitary.
The reduced row size, i.e., the number of rows of \(\Mat{\check{U}}_{\tau}\) and \(\Mat{\check{D}}_{\tau}\), is denoted by \(\check{m}_{\tau} = n_{\tau} + k\).
If no size reduction is adopted, we set \(\Mat{\Omega}_{\tau}\) to be empty and let  \(\Mat{\check{U}}_{\tau} = \Mat{U}_{\tau}\), \(\Mat{\check{D}}_{\tau} = \Mat{D}_{\tau}\) and \(\check{m}_{\tau} = m_{\tau}\).

The next is a basis elimination step.
We introduce zeros into \(\Mat{V}_{\tau} \in \complex^{n_{\tau} \times k}\) by a reverse QR factorization
\begin{equation} \label{eq:zeroing_V}
    \Mat{V}_{\tau} = \Mat{Q}_{\tau}
    \begin{bmatrix}
        \Mat{0} \\
        \Mat{\hat{V}}_{\tau; 2}
    \end{bmatrix},
\end{equation}
where \(\Mat{Q}_{\tau} \in \complex^{n_{\tau} \times n_{\tau}}\) is unitary and \(\Mat{\hat{V}}_{\tau; 2} \in \complex^{k \times k}\).
By the shared basis property of the HSS matrix, multiplying each leaf block \(\tau\) by \(\Mat{Q}_{\tau}\) on the right zeroes out corresponding columns of the node \(\tau\).
See Figure~\ref{fig:urv_process}(\subref{fig:zeroing_offdiag_cols}) for an illustration.
The diagonal block is then modified by \(\Mat{\tilde{D}}_{\tau} = \Mat{\check{D}}_{\tau} \Mat{Q}_{\tau}\).
\begin{figure}[tbhp]
    \centering
    \begin{subfigure}{0.35\textwidth}
        \begin{center}
            \begin{tikzpicture}[scale=0.45]
                \draw (0,0) rectangle (8,8);
                \draw (0,4) -- (8,4); \draw (4,0) -- (4,8);
                \foreach \j in {0, 4} {
                    \draw (\j+2, 8-\j) -- (\j+2, 4-\j);
                    \draw (8-\j,\j+2) -- (4-\j,\j+2);
                }
                \foreach \i in {0, 1, ..., 3} {
                    \pgfmathsetmacro{\ind}{\i + 4}

                    \filldraw[fill=DenseColor]
                    (2 * \i + \restdiag, 6 - 2 * \i + \restdiag)
                    rectangle
                    (2 * \i + \restdiag + \diagsz, 6 - 2 * \i + \restdiag + \diagsz);
                }
                \foreach \i in {0, 1}
                {
                    \pgfmathsetmacro{\rind}{2 * \i + 5}
                    \pgfmathsetmacro{\cind}{2 * \i + 4}

                    \filldraw[fill=DenseColor]
                    (4 * \i + \off, 4 - 4 * \i + \off)
                    rectangle
                    (4 * \i + \off + \smallsz, 4 - 4 * \i + \off + 1 + \bigsz);

                    \filldraw[fill=DenseColor]
                    (4 * \i + \off, 4 - 4 * \i + \off + 1 + \bigsz + \off)
                    rectangle
                    (4 * \i + \off + \smallsz, 4 - 4 * \i + \off + 1 + \bigsz + \off + \smallsz);

                    \filldraw[fill=white]
                    (4 * \i + \off + \smallsz + \off, 4 - 4 * \i + \off + 1 + \bigsz + \off)
                    rectangle
                    (4 * \i + \off + \smallsz + \off + 1 + \bigsz - \smallsz, 4 - 4 * \i + \off + 1 + \bigsz + \off + \smallsz);

                    \filldraw[fill=DenseColor]
                    (4 * \i + \off + \smallsz + \off + 1 + \bigsz - \smallsz , 4 - 4 * \i + \off + 1 + \bigsz + \off)
                    rectangle
                    (4 * \i + \off + \smallsz + \off + 1 + \bigsz, 4 - 4 * \i + \off + 1 + \bigsz + \off + \smallsz);
                }
                \foreach \i in {0, 1}
                {
                    \pgfmathsetmacro{\rind}{2 * \i + 4}
                    \pgfmathsetmacro{\cind}{2 * \i + 5}

                    \filldraw[fill=DenseColor]
                    (4 * \i + 2 + \off, 4 - 4 * \i + 2 + \off)
                    rectangle
                    (4 * \i + 2 + \off + \smallsz, 4 - 4 * \i + 2 + \off + 1 + \bigsz);

                    \filldraw[fill=DenseColor]
                    (4 * \i + 2 + \off, 4 - 4 * \i + 2 + \off + 1 + \bigsz + \off)
                    rectangle
                    (4 * \i + 2 + \off + \smallsz, 4 - 4 * \i + 2 + \off + 1 + \bigsz + \off + \smallsz);

                    \filldraw[fill=white]
                    (4 * \i + 2 + \off + \smallsz + \off, 4 - 4 * \i + 2 + \off + 1 + \bigsz + \off)
                    rectangle
                    (4 * \i + 2 + \off + \smallsz + \off + 1 + \bigsz - \smallsz, 4 - 4 * \i + 2 + \off + 1 + \bigsz + \off + \smallsz);

                    \filldraw[fill=DenseColor]
                    (4 * \i + 2 + \off + \smallsz + \off + 1 + \bigsz - \smallsz, 4 - 4 * \i + 2 + \off + 1 + \bigsz + \off)
                    rectangle
                    (4 * \i + 2 + \off + \smallsz + \off + 1 + \bigsz, 4 - 4 * \i + 2 + \off + 1 + \bigsz + \off + \smallsz);
                }
                {
                \pgfmathsetmacro{\rind}{2}
                \pgfmathsetmacro{\cind}{3}

                \filldraw[fill=DenseColor]
                (0 + \off + 4 * 0, 0 + \off + 4 * 0)
                rectangle
                (0 + \off + \smallsz + 4 * 0, 0 + \off + 3 + \bigsz + 4 * 0);

                \filldraw[fill=DenseColor]
                (0 + \off + 4 * 0, 0 + \off + 3 + \bigsz + \off + 4 * 0)
                rectangle
                (0 + \off + \smallsz + 4 * 0, 0 + \off + 3 + \bigsz + \off + \smallsz + 4 * 0);

                \filldraw[fill=white]
                (0 + \off + \smallsz + \off + 4 * 0, 0 + \off + 3 + \bigsz + \off + 4 * 0)
                rectangle
                (0 + \off + \smallsz + \off + 4 * 0 + \bigsz + 1 - \smallsz, 0 + \off + 3 + \bigsz + \off + \smallsz + 4 * 0);

                \filldraw[fill=DenseColor]
                (0 + \off + \smallsz + \off + 4 * 0 + \bigsz + 1 - \smallsz, 0 + \off + 3 + \bigsz + \off + 4 * 0)
                rectangle
                (0 + \off + \smallsz + \off + 4 * 0 + \bigsz + 1, 0 + \off + 3 + \bigsz + \off + \smallsz + 4 * 0);

                \filldraw[fill=white]
                (0 + \off + \smallsz + \off + 4 * 0 + \bigsz + 1, 0 + \off + 3 + \bigsz + \off + 4 * 0)
                rectangle
                (0 + \off + \smallsz + \off + 4 * 0 + 3 + \bigsz + 4 * 0 - \smallsz, 0 + \off + 3 + \bigsz + \off + \smallsz + 4 * 0);

                \filldraw[fill=DenseColor]
                (0 + \off + \smallsz + \off + 4 * 0 + 3 + \bigsz + 4 * 0 - \smallsz, 0 + \off + 3 + \bigsz + \off + 4 * 0)
                rectangle
                (0 + \off + \smallsz + \off + 3 + \bigsz + 4 * 0, 0 + \off + 3 + \bigsz + \off + \smallsz + 4 * 0);

                }

                {
                \pgfmathsetmacro{\rind}{3}
                \pgfmathsetmacro{\cind}{2}

                \filldraw[fill=DenseColor]
                (0 + \off + 4 * 1, 0 + \off + 4 * 1)
                rectangle
                (0 + \off + \smallsz + 4 * 1, 0 + \off + 3 + \bigsz + 4 * 1);

                \filldraw[fill=DenseColor]
                (0 + \off + 4 * 1, 0 + \off + 3 + \bigsz + \off + 4 * 1)
                rectangle
                (0 + \off + \smallsz + 4 * 1, 0 + \off + 3 + \bigsz + \off + \smallsz + 4 * 1);

                \filldraw[fill=white]
                (0 + \off + \smallsz + \off + 4 * 1, 0 + \off + 3 + \bigsz + \off + 4 * 1)
                rectangle
                (0 + \off + \smallsz + \off + 1 + \bigsz + 4 * 1, 0 + \off + 3 + \bigsz + \off + \smallsz + 4 * 1);

                \filldraw[fill=DenseColor]
                (0 + \off + \smallsz + \off + 1 + \bigsz + 4 * 1, 0 + \off + 3 + \bigsz + \off + 4 * 1)
                rectangle
                (0 + \off + \smallsz + \off + 1 + \bigsz + 4 * 1 + \smallsz, 0 + \off + 3 + \bigsz + \off + \smallsz + 4 * 1);

                \filldraw[fill=white]
                (0 + \off + \smallsz + \off + 1 + \bigsz + 4 * 1 + \smallsz, 0 + \off + 3 + \bigsz + \off + 4 * 1)
                rectangle
                (0 + \off + \smallsz + \off + 3 + \bigsz + 4 * 1, 0 + \off + 3 + \bigsz + \off + \smallsz + 4 * 1);

                \filldraw[fill=DenseColor]
                (0 + \off + \smallsz + \off + 3 + \bigsz + 4 * 1 - \smallsz, 0 + \off + 3 + \bigsz + \off + 4 * 1)
                rectangle
                (0 + \off + \smallsz + \off + 3 + \bigsz + 4 * 1, 0 + \off + 3 + \bigsz + \off + \smallsz + 4 * 1);}
            \end{tikzpicture}
        \end{center}
        \caption{}
        \label{fig:zeroing_offdiag_cols}
    \end{subfigure}
    \hspace{0.01\textwidth}
    \begin{subfigure}{0.35\textwidth}
        \begin{center}
            \begin{tikzpicture}[scale=0.45]
            \draw (0,0) rectangle (8,8);
            \draw (0,4) -- (8,4); \draw (4,0) -- (4,8);
            \foreach \j in {0, 4} {
                \draw (\j+2, 8-\j) -- (\j+2, 4-\j);
                \draw (8-\j,\j+2) -- (4-\j,\j+2);
            }
            \foreach \i in {0, 1, ..., 3} {
                \pgfmathsetmacro{\ind}{\i + 4}

                \filldraw[fill=white]
                (2 * \i + \restdiag, 6 - 2 * \i + \restdiag)
                rectangle
                (2 * \i + \restdiag + \diagsz - \smallsz, 6 - 2 * \i + \restdiag + \smallsz);

                \filldraw[fill=DenseColor]
                (2 * \i + \restdiag + \diagsz - \smallsz, 6 - 2 * \i + \restdiag)
                rectangle
                (2 * \i + \restdiag + \diagsz, 6 - 2 * \i + \restdiag + \smallsz);

                \filldraw[fill=DenseColor]
                (2 * \i + \restdiag, 6 - 2 * \i + \restdiag + \smallsz)
                rectangle
                (2 * \i + \restdiag + \diagsz - \smallsz, 6 - 2 * \i + \restdiag + \diagsz);

                \filldraw[fill=DenseColor]
                (2 * \i + \restdiag + \diagsz - \smallsz, 6 - 2 * \i + \restdiag + \smallsz)
                rectangle
                (2 * \i + \restdiag + \diagsz, 6 - 2 * \i + \restdiag + \diagsz);

            }
            \foreach \i in {0, 1}
            {
                \pgfmathsetmacro{\rind}{2 * \i + 5}
                \pgfmathsetmacro{\cind}{2 * \i + 4}

                \filldraw[fill=DenseColor]
                (4 * \i + \off, 4 - 4 * \i + \off)
                rectangle
                (4 * \i + \off + \smallsz, 4 - 4 * \i + \off + 1 + \bigsz);

                \filldraw[fill=DenseColor]
                (4 * \i + \off, 4 - 4 * \i + \off + 1 + \bigsz + \off)
                rectangle
                (4 * \i + \off + \smallsz, 4 - 4 * \i + \off + 1 + \bigsz + \off + \smallsz);

                \filldraw[fill=white]
                (4 * \i + \off + \smallsz + \off, 4 - 4 * \i + \off + 1 + \bigsz + \off)
                rectangle
                (4 * \i + \off + \smallsz + \off + 1 + \bigsz - \smallsz, 4 - 4 * \i + \off + 1 + \bigsz + \off + \smallsz);

                \filldraw[fill=DenseColor]
                (4 * \i + \off + \smallsz + \off + 1 + \bigsz - \smallsz , 4 - 4 * \i + \off + 1 + \bigsz + \off)
                rectangle
                (4 * \i + \off + \smallsz + \off + 1 + \bigsz, 4 - 4 * \i + \off + 1 + \bigsz + \off + \smallsz);
            }
            \foreach \i in {0, 1}
            {
                \pgfmathsetmacro{\rind}{2 * \i + 4}
                \pgfmathsetmacro{\cind}{2 * \i + 5}

                \filldraw[fill=DenseColor]
                (4 * \i + 2 + \off, 4 - 4 * \i + 2 + \off)
                rectangle
                (4 * \i + 2 + \off + \smallsz, 4 - 4 * \i + 2 + \off + 1 + \bigsz);

                \filldraw[fill=DenseColor]
                (4 * \i + 2 + \off, 4 - 4 * \i + 2 + \off + 1 + \bigsz + \off)
                rectangle
                (4 * \i + 2 + \off + \smallsz, 4 - 4 * \i + 2 + \off + 1 + \bigsz + \off + \smallsz);

                \filldraw[fill=white]
                (4 * \i + 2 + \off + \smallsz + \off, 4 - 4 * \i + 2 + \off + 1 + \bigsz + \off)
                rectangle
                (4 * \i + 2 + \off + \smallsz + \off + 1 + \bigsz - \smallsz, 4 - 4 * \i + 2 + \off + 1 + \bigsz + \off + \smallsz);

                \filldraw[fill=DenseColor]
                (4 * \i + 2 + \off + \smallsz + \off + 1 + \bigsz - \smallsz, 4 - 4 * \i + 2 + \off + 1 + \bigsz + \off)
                rectangle
                (4 * \i + 2 + \off + \smallsz + \off + 1 + \bigsz, 4 - 4 * \i + 2 + \off + 1 + \bigsz + \off + \smallsz);
            }
            {
            \pgfmathsetmacro{\rind}{2}
            \pgfmathsetmacro{\cind}{3}

            \filldraw[fill=DenseColor]
            (0 + \off + 4 * 0, 0 + \off + 4 * 0)
            rectangle
            (0 + \off + \smallsz + 4 * 0, 0 + \off + 3 + \bigsz + 4 * 0);

            \filldraw[fill=DenseColor]
            (0 + \off + 4 * 0, 0 + \off + 3 + \bigsz + \off + 4 * 0)
            rectangle
            (0 + \off + \smallsz + 4 * 0, 0 + \off + 3 + \bigsz + \off + \smallsz + 4 * 0);

            \filldraw[fill=white]
            (0 + \off + \smallsz + \off + 4 * 0, 0 + \off + 3 + \bigsz + \off + 4 * 0)
            rectangle
            (0 + \off + \smallsz + \off + 4 * 0 + \bigsz + 1 - \smallsz, 0 + \off + 3 + \bigsz + \off + \smallsz + 4 * 0);

            \filldraw[fill=DenseColor]
            (0 + \off + \smallsz + \off + 4 * 0 + \bigsz + 1 - \smallsz, 0 + \off + 3 + \bigsz + \off + 4 * 0)
            rectangle
            (0 + \off + \smallsz + \off + 4 * 0 + \bigsz + 1, 0 + \off + 3 + \bigsz + \off + \smallsz + 4 * 0);

            \filldraw[fill=white]
            (0 + \off + \smallsz + \off + 4 * 0 + \bigsz + 1, 0 + \off + 3 + \bigsz + \off + 4 * 0)
            rectangle
            (0 + \off + \smallsz + \off + 4 * 0 + 3 + \bigsz + 4 * 0 - \smallsz, 0 + \off + 3 + \bigsz + \off + \smallsz + 4 * 0);

            \filldraw[fill=DenseColor]
            (0 + \off + \smallsz + \off + 4 * 0 + 3 + \bigsz + 4 * 0 - \smallsz, 0 + \off + 3 + \bigsz + \off + 4 * 0)
            rectangle
            (0 + \off + \smallsz + \off + 3 + \bigsz + 4 * 0, 0 + \off + 3 + \bigsz + \off + \smallsz + 4 * 0);

            }

            {
            \pgfmathsetmacro{\rind}{3}
            \pgfmathsetmacro{\cind}{2}

            \filldraw[fill=DenseColor]
            (0 + \off + 4 * 1, 0 + \off + 4 * 1)
            rectangle
            (0 + \off + \smallsz + 4 * 1, 0 + \off + 3 + \bigsz + 4 * 1);

            \filldraw[fill=DenseColor]
            (0 + \off + 4 * 1, 0 + \off + 3 + \bigsz + \off + 4 * 1)
            rectangle
            (0 + \off + \smallsz + 4 * 1, 0 + \off + 3 + \bigsz + \off + \smallsz + 4 * 1);

            \filldraw[fill=white]
            (0 + \off + \smallsz + \off + 4 * 1, 0 + \off + 3 + \bigsz + \off + 4 * 1)
            rectangle
            (0 + \off + \smallsz + \off + 1 + \bigsz + 4 * 1, 0 + \off + 3 + \bigsz + \off + \smallsz + 4 * 1);

            \filldraw[fill=DenseColor]
            (0 + \off + \smallsz + \off + 1 + \bigsz + 4 * 1, 0 + \off + 3 + \bigsz + \off + 4 * 1)
            rectangle
            (0 + \off + \smallsz + \off + 1 + \bigsz + 4 * 1 + \smallsz, 0 + \off + 3 + \bigsz + \off + \smallsz + 4 * 1);

            \filldraw[fill=white]
            (0 + \off + \smallsz + \off + 1 + \bigsz + 4 * 1 + \smallsz, 0 + \off + 3 + \bigsz + \off + 4 * 1)
            rectangle
            (0 + \off + \smallsz + \off + 3 + \bigsz + 4 * 1, 0 + \off + 3 + \bigsz + \off + \smallsz + 4 * 1);

            \filldraw[fill=DenseColor]
            (0 + \off + \smallsz + \off + 3 + \bigsz + 4 * 1 - \smallsz, 0 + \off + 3 + \bigsz + \off + 4 * 1)
            rectangle
            (0 + \off + \smallsz + \off + 3 + \bigsz + 4 * 1, 0 + \off + 3 + \bigsz + \off + \smallsz + 4 * 1);
            }
            \end{tikzpicture}
        \end{center}
        \caption{}
        \label{fig:QR_diag_blocks}
    \end{subfigure}
    \caption{An illustration of the URV factorization. Left: Introducing zeros into off-diagonal columns. Right: QR factorization on diagonal blocks. }
    \label{fig:urv_process}
\end{figure}

By partitioning the diagonal block as \(\Mat{\tilde{D}}_{\tau} =
    \begin{bmatrix}
        \Mat{\tilde{D}}_{\tau; 1, 1} & \Mat{\tilde{D}}_{\tau; 1, 2} \\
        \Mat{\tilde{D}}_{\tau; 2, 1} & \Mat{\tilde{D}}_{\tau; 2, 2} \\
    \end{bmatrix}\),
our goal is to eliminate the block \(\Mat{\tilde{D}}_{\tau; 2, 1}\) by a unitary transformation.
This can be done by a QR factorization of the first column block:
\begin{equation*}
    \begin{bmatrix}
        \Mat{\tilde{D}}_{\tau; 1, 1} \\
        \Mat{\tilde{D}}_{\tau; 2, 1} \\
    \end{bmatrix}
    = \Mat{P}_{\tau}
    \begin{bmatrix}
        \Mat{\hat{D}}_{\tau; 1, 1} \\
        \Mat{0} \\
    \end{bmatrix},
\end{equation*}
where \(\Mat{P}_{\tau} \in \complex^{\check{m}_{\tau} \times \check{m}_{\tau}}\) is unitary.
After multiplying \(\Mat{P}_{\tau}^{\herm}\) on the left, the diagonal block transforms into
\begin{equation*}
    \Mat{\hat{D}}_{\tau} = \Mat{P}_{\tau}^{\herm} \Mat{\tilde{D}}_{\tau} =
    \begin{bmatrix}
        \Mat{\hat{D}}_{\tau; 1, 1} & \Mat{\hat{D}}_{\tau; 1, 2} \\
        \Mat{0} & \Mat{\hat{D}}_{\tau; 2, 2}
    \end{bmatrix}
\end{equation*}
and the basis matrix \(\Mat{\check{U}}_{\tau}\) is modified to
\begin{equation*}
    \Mat{P}_{\tau}^{\herm} \Mat{\check{U}}_{\tau} =
    \begin{bmatrix}
        \Mat{\hat{U}}_{\tau; 1} \\
        \Mat{\hat{U}}_{\tau; 2}
    \end{bmatrix}.
\end{equation*}
Figure~\ref{fig:urv_process}(\subref{fig:QR_diag_blocks}) illustrates the results after the elimination on diagonal blocks.

Now, it is remarkable that if we partition \(\tau = \tau_1 \sqcup \tau_2\) conformably for all leaf nodes \(\tau\), the subproblem with respect to \(\tau_2\) is entirely decoupled from \(\tau_1\).
Moreover, for a nonleaf node \(\tau\) with children \(\alpha_{1}\) and \(\alpha_{2}\), by merging appropriate blocks as
\begin{equation} \label{eq:urv_merge}
    \begin{gathered}
        \Mat{D}_{\tau}^{+} =
        \begin{bmatrix}
            \Mat{\hat{D}}_{\alpha_{1}; 2, 2} & \Mat{\hat{U}}_{\alpha_{1}; 2} \Mat{B}_{\alpha_{1}, \alpha_{2}} \Mat{\hat{V}}_{\alpha_{2}; 2}^{\herm} \\
            \Mat{\hat{U}}_{\alpha_{2}; 2} \Mat{B}_{\alpha_{2}, \alpha_{1}} \Mat{\hat{V}}_{\alpha_{1}; 2}^{\herm} & \Mat{\hat{D}}_{\alpha_{2}; 2, 2} \\
        \end{bmatrix}, \\
        \Mat{U}_{\tau}^{+} =
        \begin{bmatrix}
            \Mat{\hat{U}}_{\alpha_{1}; 2} & \\
            & \Mat{\hat{U}}_{\alpha_{2}; 2}  \\
        \end{bmatrix}
        \Mat{U}_{\tau}
        , \quad
        \Mat{V}_{\tau}^{+} =
        \begin{bmatrix}
            \Mat{\hat{V}}_{\alpha_{1}; 2} & \\
            & \Mat{\hat{V}}_{\alpha_{2}; 2} \\
        \end{bmatrix}
        \Mat{V}_{\tau},
    \end{gathered}
\end{equation}
then the children \(\alpha_{1}\) and \(\alpha_{2}\) can be ``eliminated'' and \(\tau\) becomes a new leaf with the associated generators defined in~\eqref{eq:urv_merge}.
After that, the remaining matrix still has the HSS structure but with a maximum level \(L - 1\).
This enables us to recursively factor the merged HSS matrix~(with \(\Mat{D}_{\tau}\), \(\Mat{U}_{\tau}\), \(\Mat{V}_{\tau}\) replaced by \(\Mat{D}_{\tau}^{+}\), \(\Mat{U}_{\tau}^{+}\), \(\Mat{V}_{\tau}^{+}\)) until the root node is the only remaining node.
At that point, we only need to perform a QR decomposition \(\Mat{D}_{\rho}^{+} = \Mat{P}_{\rho} \Mat{\hat{D}}_{\rho}\).

Once the URV factorization is computed, the solution process can be divided into two steps similar to the application of QR factorization.
The first step applies unitary transforms to the right-hand side, and the second step uses backward substitution to solve a sequence of hierarchical triangular systems.

The first step involves a bottom-up traversal of the HSS tree.
Starting from each leaf node \(\tau\), let \(\Mat{f}_{\tau}\) denote the corresponding block of the right-hand side.
If size reduction is applied to \(\tau\), we define
\(\Mat{\check{f}}_{\tau} = \Mat{\Omega}_{\tau}^{\herm} \Mat{f}_{\tau} =
    \begin{bmatrix}
        \Mat{\check{f}}_{\tau; 1} \\
        \Mat{\check{f}}_{\tau; 2} \\
    \end{bmatrix}\) and
\(\Mat{\hat{f}}_{\tau} = \Mat{P}_{\tau}^{\herm} \Mat{\check{f}}_{\tau; 1} =
    \begin{bmatrix}
        \Mat{\hat{f}}_{\tau; 1} \\
        \Mat{\hat{f}}_{\tau; 2} \\
    \end{bmatrix}\).
If there is no size reduction step, the first equation is skipped and \(\Mat{\check{f}}_{\tau; 1}\) is replaced by \(\Mat{f}_{\tau}\) in the second equation.
Next, we merge the vectors on each nonleaf node \(\tau\) with children \(\alpha_{1}\) and \(\alpha_{2}\) by
\(\Mat{f}_{\tau}^{+} =
    \begin{bmatrix}
        \Mat{\hat{f}}_{\alpha_{1}; 2} \\
        \Mat{\hat{f}}_{\alpha_{2}; 2} \\
    \end{bmatrix}\).
We continue applying the same procedure~(with \(\Mat{f}_{\tau}\) replaced by \(\Mat{f}_{\tau}^{+}\)) until \(\tau\) is the root node.

The second step is a top-down procedure.
Beginning at the root \(\rho\), we solve the triangular system \(\Mat{\hat{D}}_{\rho} \Mat{c}_{\rho} = \Mat{\hat{f}}_{\rho}\).
For its children \(\alpha_{1}\) and \(\alpha_{2}\), we partition the solution as
\(\Mat{c}_{\rho} =
    \begin{bmatrix}
        \Mat{\hat{c}}_{\alpha_{1}; 2} \\
        \Mat{\hat{c}}_{\alpha_{2}; 2} \\
    \end{bmatrix}\).
We then compute \(\Mat{\hat{g}}_{\alpha_{1}} = \Mat{B}_{\alpha_{1}, \alpha_{2}} \Mat{\hat{V}}_{\alpha_{2}; 2}^{\herm} \Mat{\hat{c}}_{\alpha_{2}; 2}\) and \(\Mat{\hat{g}}_{\alpha_{2}} = \Mat{B}_{\alpha_{2}, \alpha_{1}} \Mat{\hat{V}}_{\alpha_{1}; 2}^{\herm} \Mat{\hat{c}}_{\alpha_{1}; 2}\).
For a nonroot node \(\tau\), we obtain \(\Mat{\hat{c}}_{\tau; 1}\) by solving the triangular system \(\Mat{\hat{D}}_{\tau; 1, 1} \Mat{\hat{c}}_{\tau; 1} = \Mat{\hat{f}}_{\tau; 1} - \Mat{\hat{D}}_{\tau; 1, 2} \Mat{\hat{c}}_{\tau; 2} - \Mat{\hat{U}}_{\tau; 1} \Mat{\hat{g}}_{\tau}\) and update
\(\Mat{c}_{\tau} = \Mat{Q}_{\tau}
    \begin{bmatrix}
        \Mat{\hat{c}}_{\tau; 1} \\
        \Mat{\hat{c}}_{\tau; 2} \\
    \end{bmatrix}\).
If \(\tau\) is a nonleaf node with children \(\alpha_{1}\) and \(\alpha_{2}\), we partition \(\Mat{c}_{\tau}\) conformably as \(\Mat{c}_{\tau} =
\begin{bmatrix}
    \Mat{\hat{c}}_{\alpha_{1}; 2} \\
    \Mat{\hat{c}}_{\alpha_{2}; 2} \\
\end{bmatrix}\)
and compute \(\Mat{\hat{g}}_{\alpha_{1}} = \Mat{B}_{\alpha_{1}, \alpha_{2}} \Mat{\hat{V}}_{\alpha_{2}; 2}^{\herm} \Mat{\hat{c}}_{\alpha_{2}; 2} + \Mat{R}_{\alpha_{1}} \Mat{\hat{g}}_{\tau}\) and \(\Mat{\hat{g}}_{\alpha_{2}} = \Mat{B}_{\alpha_{2}, \alpha_{1}} \Mat{\hat{V}}_{\alpha_{1}; 2}^{\herm} \Mat{\hat{c}}_{\alpha_{1}; 2} + \Mat{R}_{\alpha_{2}} \Mat{\hat{g}}_{\tau}\),
where \(\Mat{U}_{\tau} =
    \begin{bmatrix}
        \Mat{R}_{\alpha_{1}} \\
        \Mat{R}_{\alpha_{2}} \\
    \end{bmatrix}\).
Repeating this process until we reach the leaf nodes, the solution \(\Mat{c}\) is obtained by collecting vectors \(\Mat{c}_{\tau}\) for all leaf nodes \(\tau\).

\subsection{A superfast direct solver of type-II INUDFT problem} \label{subsec:inudft2}

The type-II NUDFT matrix \(\Mat{A}\), defined by \(\Mat{A}(j, k) = \e^{-2 \pi \imag x_{j} k} = \gamma_{j}^{k}\) with \(\gamma_{j} = \e^{-2 \pi \imag x_{j}}\) on the unit circle \(\unitCircle\), is a Vandermonde matrix and satisfies the Sylvester equation
\begin{equation} \label{eq:nudft2_Sylvester_A}
    \Mat{\Gamma} \Mat{A} - \Mat{A} \Mat{C} = \Mat{u} \Mat{e}_{N - 1}^{\top},
\end{equation}
where \(\Mat{\Gamma} = \diag(\gamma_0, \gamma_1, \dotsc, \gamma_{M - 1}) \in \complex^{M \times M}\) is diagonal,
\(\Mat{C} = \begin{bmatrix}
\Mat{0} & 1 \\
\Mat{I}_{N - 1} & \Mat{0}\\
\end{bmatrix} \in \complex^{N \times N}\) is the shift-down circulant matrix, and \(u_{j} = \gamma_{j}^N - 1\).
Here \(\Mat{e}_{N - 1}\) is the canonical unit vector in \(\complex^{N}\).
The circulant matrix \(\Mat{C}\) can be diagonalized using the DFT matrix \(\Mat{F}\), expressed as \(\Mat{C} = \Mat{F}^{-1} \Mat{D} \Mat{F}\), where \(\Mat{D} = \diag(\xi_{0}, \xi_{1}, \dotsc, \xi_{N - 1})\) is a diagonal matrix with \(\xi_{k} = \zeta_{N}^{k}\) for \(k = 0, 1, \dotsc, N - 1\).
Here \(\zeta_{N} = \e^{-2 \pi \imag / N}\) is the \(N\)-th root of unity.
By multiplying \(\Mat{F}^{-1}\) to the right of ~\eqref{eq:nudft2_Sylvester_A} and denoting \(\Mat{\tilde{A}} = \Mat{A} \Mat{F}^{-1}\), the equation becomes
\begin{equation} \label{eq:nudft2_Sylvester_tildeA}
    \Mat{\Gamma} \Mat{\tilde{A}} - \Mat{\tilde{A}} \Mat{D} = \Mat{u} \Mat{v}^{\top},
\end{equation}
where the entries of \(\Mat{v}\) are given by \(v_{k} = \zeta_{N}^{k} / N = \xi_{k} / N\) for \(k = 0, 1, \dotsc, N - 1\).

In~\cite{Wilber_Epperly_Barnett_2025}, the authors utilize the displacement structure~\eqref{eq:nudft2_Sylvester_tildeA} to establish the low-rank property of \(\Mat{\tilde{A}}\) and provide an estimate of the numerical rank~\cite{Beckermann_Townsend_2017}.
Specifically, let the \(\varepsilon\)-rank of a matrix \(\Mat{M}\) be defined as \(\rank_{\varepsilon}(\Mat{M}) \coloneq \min\{k: \|\Mat{X} - \Mat{M}\|_2 \leq \varepsilon \|\Mat{M}\|_2, \rank(\Mat{X}) = k\}\).
It is shown that the \(\varepsilon\)-rank of the HSS blocks of \(\Mat{\tilde{A}}\) admits an upper bound of \(k = \bigO\bigl(\log(1 / \varepsilon) \log N\bigr)\) (Theorem 3.2 in~\cite{Wilber_Epperly_Barnett_2025}).
See~\cite{Beckermann_Kressner_Wilber_2025} for a more general result for the Toeplitz-like matrices.
Additionally, by exploiting the displacement structure, the HSS matrix approximation can be constructed using the \textit{factorized alternating direction implicit}~(fADI) method~\cite{Benner_Li_Truhar_2009, Wilber_2021}, resulting in an overall complexity of \(\bigO\bigl(k^{2} (M + N)\bigr)\).

Here we present an alternative perspective on the low-rank property.
Direct calculations show that the entries of \(\Mat{\tilde{A}}\) are explicitly given by
\begin{equation} \label{eq:nudft2_kernel_matrix}
    \Mat{\tilde{A}}(j, k) = \frac{u_{j} v_{k}}{\gamma_{j} - \xi_{k}} = \frac{1}{N} \frac{(\gamma_{j}^{N} - 1) \xi_{k}}{\gamma_{j} - \xi_{k}}.
\end{equation}
Consequently, \(\Mat{\tilde{A}}\) is a kernel matrix corresponding to kernel function
\begin{equation*}
    K(z, w) = \frac{1}{N} \frac{(z^N - 1) w}{z - w}
\end{equation*}
evaluated at the points \(\{\gamma_{j}\}\) and \(\{\xi_{k}\}\) for rows and columns.
Using Taylor expansion or pole expansion, one can show that for well-separated \(z\) and \(w\), the \(\varepsilon\)-rank of the kernel matrix is bounded by \(\bigO(\log(1 / \varepsilon))\), independent of size~\cite{Hackbusch_Khoromskij_2000}.
Then using the technique in~\cite{Ho_Greengard_2012}, by recursively applying the well-separated block rank bound to the HSS blocks, we can show that the \(\varepsilon\)-rank of the HSS blocks is bounded by \(\bigO(\log(1 / \varepsilon) \log N)\).
See the discussion in section~4 in~\cite{Ho_Greengard_2012} for more details.
Consequently, \(\Mat{\tilde{A}}\) can be approximated by an HSS matrix whose HSS rank is bounded by \(k = \bigO(\log(1 / \varepsilon) \log N)\).
This is consistent with the discussion in~\cite{Wilber_Epperly_Barnett_2025}.

To conclude, the type-II NUDFT matrix \(\Mat{A}\) is approximated by \(\Mat{A}_{\fast} \coloneq \Mat{\tilde{A}}_{\HSS} \Mat{F}\), where \(\Mat{\tilde{A}}_{\HSS}\) is the HSS approximation of \(\Mat{\tilde{A}}\).
The type-II direct solver is then obtained by applying the least-squares solver in section~\ref{subsec:ls_solver_hss} and the iFFT.
We summarize the whole procedure in Algorithm~\ref{alg:inudft2}.
The algorithm consists of three stages: construction, factorization, and solution.
In the construction and factorization stages, the HSS compression of \(\Mat{\tilde{A}}\) is constructed by the fADI-based method and factorized by the URV factorization, respectively.
Once the factorization is finished, the solution \(\Mat{c}\) can be obtained using the solution process of the URV factorization and the iFFT.
The complexity of the construction and factorization is \(\bigO\bigl(k^{2} (M + N)\bigr)\), while the solution stage has a complexity of \(\bigO\bigl(k M + N \log N\bigr)\), where \(k\) is the HSS rank of \(\Mat{\tilde{A}}_{\HSS}\).
Using the bound \(k = \bigO(\log N)\) and \(M \geq N\), the type-II solver has the complexity of \(\bigO\bigl((M + N) \log^{2} N\bigr) = \bigO\bigl(M \log^{2} N\bigr)\)  for construction and factorization and \(\bigO\bigl((M + N) \log N\bigr) = \bigO(M \log N)\) for solution.
\begin{algorithm}[htbp]
  \caption{A direct solver for type-II INUDFT problem. \label{alg:inudft2}}
  \begin{algorithmic}[1]
    \REQUIRE{Sample points \(\{x_{j}\}_{j = 0}^{M - 1}\), number of frequencies \(N\), target values \(\{f_j\}_{j = 0}^{M - 1}\), accuracy parameter \(\varepsilon \in (0, 1)\).}
    \ENSURE{Coefficients \(\{c_{k}\}_{k = 0}^{N - 1}\).}
    \STATE{Approximate \(\Mat{\tilde{A}} = \Mat{A} \Mat{F}^{-1}\) by an HSS matrix \(\Mat{\tilde{A}}_{\HSS}\) using the fADI-based construction with accuracy parameter \(\varepsilon\).}
    \STATE{Compute the URV factorization of \(\Mat{\tilde{A}}_{\HSS}\).}
    \STATE{Solve \(\Mat{v} = \arg \min_{\Mat{w}} \|\Mat{\tilde{A}}_{\HSS} \Mat{w} - \Mat{f}\|_2\) by the HSS least-squares solver.}
    \STATE{Compute \(\Mat{c} = \Mat{F}^{-1} \Mat{v}\) using the iFFT.}
  \end{algorithmic}
\end{algorithm}

\section{A superfast direct solver of the type-III INUDFT problem} \label{sec:inudft3}

In this section, we propose a superfast direct solver for the type-III INUDFT problem.
We first give some intuitions and motivations in section~\ref{subsec:inudft3_motivations}, and then present the detailed algorithm and its implementation in section~\ref{subsec:inudft3_algorithm}.

\subsection{Motivations} \label{subsec:inudft3_motivations}

Our motivation comes from the expansion of the nonharmonic Fourier basis in terms of the Fourier basis.
The square integrable functions on \([0, 1]\) form a Hilbert space \(L^{2}[0, 1]\) with the inner product defined by \(\langle f, g \rangle \coloneq \int_{0}^{1} f(x) \overline{g(x)} \diff x\).
A typically used orthonormal basis of \(L^{2}[0, 1]\) is given by \(\{\e^{-2 \pi \imag \ell x}\}_{\ell \in \integer}\).
We begin by expanding the function \(\e^{-2 \pi \imag \omega_{k} x}\) in terms of this basis:
\begin{equation} \label{eq:nonuniform_l2_expansion}
    \e^{-2 \pi \imag \omega_{k} x} = \sum_{\ell \in \integer} \langle \e^{-2 \pi \imag \omega_{k} x}, \e^{-2 \pi \imag \ell x} \rangle \e^{-2 \pi \imag \ell x} = \sum_{\ell \in \integer} Q_{\ell, k} \e^{-2 \pi \imag \ell x},
\end{equation}
where the coefficients are given by
\begin{equation} \label{eq:Q_entry}
    Q_{\ell, k} = \int_{0}^{1} \e^{2 \pi \imag (\ell - \omega_{k}) x} \diff x
    = \e^{\pi \imag (\ell- \omega_{k})} \frac{\sin \bigl(\pi (\ell- \omega_{k})\bigr)}{\pi (\ell- \omega_{k})} =  G(\ell- \omega_{k}).
\end{equation}
The assembled matrix \(\Mat{Q}\) can be viewed as a kernel matrix
with the kernel function defined by
\begin{equation} \label{eq:nudft_kernel}
    G(\omega) \coloneq \e^{\pi \imag \omega} \frac{\sin(\pi \omega)}{\pi \omega}.
\end{equation}
Since the kernel \(G\) is smooth and not highly oscillatory, the matrix \(\Mat{Q}\) is expected to have a low-rank property for its off-diagonal blocks.
Suppose \(N\) is even and \(\sigma \geq 1\) is some integer, and consider the index set \(\indexset{I}_{N; \sigma} = \{\ell \in \integer: -(\sigma N - N) / 2 \leq \ell \leq (\sigma N + N) / 2 - 1\}\), which is an expansion of \(\{0, 1, \dotsc, N - 1\}\) about \(N / 2\) by a factor of \(\sigma\).
Truncating the series in~\eqref{eq:nonuniform_l2_expansion} leads to
\begin{equation*}
    \e^{-2 \pi \imag \omega_{k} x} = \sum_{\ell \in \integer} Q_{\ell, k} \e^{-2 \pi \imag \ell x} = \sum_{\ell \in \indexset{I}_{N; \sigma}} Q_{\ell, k} \e^{-2 \pi \imag \ell x} + \sum_{\ell \notin \indexset{I}_{N; \sigma}} Q_{\ell, k} \e^{-2 \pi \imag \ell x},
\end{equation*}
where the first and second terms in the right-hand side serve as the approximation and the error, respectively.
The \(L^{2}\) norm of the error, denoted as \(\err_{k}(x)\), is given by Parseval's identity:
\begin{equation*}
    \|\err_{k}(\cdot)\|_{L^{2}}^2 = \sum_{\ell \notin \indexset{I}_{N; \sigma}} |Q_{\ell, k}|^2 = \sum_{\ell \notin \indexset{I}_{N; \sigma}} \frac{\sin^2 \bigl(\pi (\ell - \omega_{k})\bigr)}{\pi^2 (\ell - \omega_{k})^2}.
\end{equation*}
Assume \(\{\omega_{k} = k + \delta_{k}\}\) is a perturbation of the continuous integers, where \(|\delta_{k}| \leq \alpha\) for some constant \(\alpha \in [0, 1 / 2)\).
Then the error can be estimated as
\begin{equation*}
    \begin{aligned}
        \|\err_{k}(\cdot)\|_{L^{2}}^2
        & = \sum_{\ell = -\infty}^{-(\sigma N - N) / 2 - 1} \frac{\sin^2 \bigl(\pi (\ell - \omega_{k})\bigr)}{\pi^2 (\ell - \omega_{k})^2}
        + \sum_{\ell = (\sigma N + N) / 2}^{\infty} \frac{\sin^2 \bigl(\pi (\ell - \omega_{k})\bigr)}{\pi^2 (\ell - \omega_{k})^2} \\
        & = \sum_{\ell = -\infty}^{-(\sigma N - N) / 2 - 1} \frac{\sin^2 (\pi \delta_{k})}{\pi^2 (\ell - \omega_{k})^2}
        + \sum_{\ell = (\sigma N + N) / 2}^{\infty} \frac{\sin^2 (\pi \delta_{k})}{\pi^2 (\ell - \omega_{k})^2} \\
        & \leq 2 \sum_{s = (\sigma - 1) N / 2}^{\infty} \frac{\sin^2 (\pi \alpha)}{\pi^2 (s + 1 - \alpha)^2}. \\
    \end{aligned}
\end{equation*}
By the integral estimate of the tail and noting that the integrand is a convex function, we have
\begin{equation*}
    \begin{aligned}
        \sum_{s = S}^{\infty} \frac{1}{(s + 1 - \alpha)^2}
         & \leq \sum_{s = S}^{\infty} \int_{s - 1 / 2}^{s + 1 / 2} \frac{1}{(x + 1 - \alpha)^2} \diff x \\
         & \leq \int_{S - 1 / 2}^{\infty} \frac{1}{(x + 1 - \alpha)^2} \diff x
        = \frac{1}{S + 1 / 2 - \alpha}.
    \end{aligned}
\end{equation*}
Therefore,
\begin{equation*}
    \|\err_{k}(\cdot)\|_{L^{2}}^2 \leq \frac{2 \sin^2 (\pi \alpha)}{\pi^2 ((\sigma - 1) N / 2 + 1 / 2 - \alpha)}.
\end{equation*}

For our type-III NUDFT matrix, each entry can be approximated by
\begin{equation*}
    \Mat{A}(j, k) = \e^{-2 \pi \imag x_{j} \omega_{k}}
    = \sum_{\ell \in \indexset{I}_{N; \sigma}} \e^{-2 \pi \imag x_{j} \ell} Q_{\ell, k}
    + \sum_{\ell \notin \indexset{I}_{N; \sigma}} \e^{-2 \pi \imag x_{j} \ell} Q_{\ell, k}.
\end{equation*}
The approximation can be reformulated in matrix form as
\begin{equation} \label{eq:nudft3_matrix_approx}
    \Mat{A} = \Mat{B} \Mat{Q} + \Mat{E},
\end{equation}
where \(\Mat{B} = (\e^{-2 \pi \imag x_{j} \ell})_{j, \ell}\), \(\Mat{Q} = (Q_{\ell, k})_{\ell, k}\), and \(\Mat{E} = (\sum_{\ell \notin \indexset{I}_{N; \sigma}} Q_{\ell, k} \e^{-2 \pi \imag \ell x_{j}})_{j, k}\) are the corresponding type-II NUDFT matrix, coefficient matrix, and error matrix, respectively.
To estimate the approximation error, we assume that the sample points \(\{x_{j}\}\) are i.i.d.\ uniform random variables in \([0, 1)\).
The statistical properties of the random matrix \(\Mat{E}\) can be obtained through direct calculations.
Precisely, we have \(\expectation \Mat{E}(j, k) = \sum_{\ell \notin \indexset{I}_{N; \sigma}} Q_{\ell, k} \int_{0}^{1} \e^{-2 \pi \imag \ell x} \diff x = 0\) and
\begin{equation*}
    \begin{aligned}
        \expectation \|\Mat{E}\|_{\fro}^2
        & = \expectation \sum_{j = 0}^{M - 1} \sum_{k = 0}^{N - 1} \biggl| \sum_{\ell \notin \indexset{I}_{N; \sigma}} Q_{\ell, k} \e^{-2 \pi \imag \ell x_{j}} \biggr|^2
        = M \int_{0}^{1} \sum_{k = 0}^{N - 1} \biggl| \sum_{\ell \notin \indexset{I}_{N; \sigma}} Q_{\ell, k} \e^{2 \pi \imag \ell x} \biggr|^2 \diff x \\
        & = M \sum_{k = 0}^{N - 1} \sum_{\ell \notin \indexset{I}_{N; \sigma}} \bigl|Q_{\ell, k}\bigr|^2
        \leq M \sum_{k = 0}^{N - 1} \|\err_{k}(\cdot)\|_{L^{2}}^2 \\
        & \leq M N \frac{2 \sin^2 (\pi \alpha)}{\pi^2 ((\sigma - 1) N / 2 + 1 / 2 - \alpha)} = \frac{2 \sin^2 (\pi \alpha)}{\pi^2 ((\sigma - 1) N / 2 + 1 / 2 - \alpha)} \|\Mat{A}\|_{\fro}^2.
    \end{aligned}
\end{equation*}
We summarize the results in the following proposition by noting that \(\|\Mat{A}\|_{\fro}^2 = M N\).
\begin{proposition} \label{prop:nudft3_matrix_approx_error}
    Suppose \(\Mat{A}\in \complex^{M \times N}\) is the type-III NUDFT matrix defined as \(\Mat{A}(j, k) = \e^{-2 \pi \imag x_{j} \omega_{k}}\), where \(x_{j}\) are i.i.d.\ uniform random variables in \([0, 1)\) and \(\omega_{k} = k + \delta_{k}\) are perturbations of the continuous integers with \(|\delta_{k}| \leq \alpha\) for some constant \(\alpha \in [0, 1 / 2)\).
    Let \(\sigma \geq 1\), and let \(\Mat{B}\) and \(\Mat{Q}\) be the matrices defined in~\eqref{eq:nudft3_matrix_approx}.
    Then the error matrix \(\Mat{E}\) satisfies \(\expectation \Mat{E} = \Mat{0}\) and \(\expectation \|\Mat{E}\|_{\fro}^2 \leq \frac{2 \sin^2 (\pi \alpha)}{\pi^2 ((\sigma - 1) N / 2 + 1 / 2 - \alpha)} \|\Mat{A}\|_{\fro}^2\).
\end{proposition}

\subsection{Algorithm} \label{subsec:inudft3_algorithm}

While the approximation derived from~\eqref{eq:nudft3_matrix_approx} is straightforward, it is not the form used in our practical algorithm.
One main reason is that the error of this approximation can be quite large.
Our strategy is to replace the matrix \(\Mat{Q}\) by the best candidate in the approximation \(\|\Mat{A} - \Mat{B} \Mat{H}\|_{\fro}\), i.e., \(\Mat{H} = \Mat{B}^{\pinv} \Mat{A}\).
We will later demonstrate that this leads to a reprojection of the error matrix \(\Mat{E}\) against the type-II NUDFT matrix \(\Mat{B}\).
In this section, we describe our practical superfast direct solver for the type-III INUDFT problem.
The overall algorithm consists of four stages: preliminary, construction, factorization, and solution.
We will discuss each stage in detail in the rest of this section, and the overall pseudocode can be found in Algorithm~\ref{alg:inudft3}.
\begin{algorithm}[htbp]
  \caption{A direct solver for type-III INUDFT problem. \label{alg:inudft3}}
  \begin{algorithmic}[1]
    \REQUIRE{Sample points \(\{x_{j}\}_{j = 0}^{M - 1}\), frequencies \(\{\omega_{k}\}_{k = 0}^{N - 1}\), target values \(\{f_j\}_{j = 0}^{M - 1}\), accuracy parameter \(\varepsilon \in (0, 1)\), and rank parameter \(k\).}
    \ENSURE{Coefficients \(\{c_{k}\}_{k = 0}^{N - 1}\).}
    \STATE{Construct a type-II direct solver \(\Mat{B} \approx \Mat{B}_{\fast} = \Mat{\tilde{B}}_{\HSS} \Mat{F}\) by Algorithm~\ref{alg:inudft2} for the type-II NUDFT matrix \(\Mat{B}(j, k) = \e^{-2 \pi \imag x_{j} k}\) with accuracy parameter \(\varepsilon\).}
    \STATE{Set \(r = k + p\)~(\(p = 5\) or \(p = 10\)) and \(s \geq \max(r + m, 3r)\), where \(m\) is the max leaf size in the HSS tree corresponding to \(\Mat{H}_{\HSS}\).}
    \STATE{Generate \(N \times s\) Gaussian random matrices \(\Mat{\Omega}\) and \(\Mat{\Psi}\).}
    \STATE{Compute \(\Mat{Y} = \Mat{B}_{\fast}^{\pinv} \Mat{A} \Mat{\Omega}\) and \(\Mat{Z} = \Mat{A}^{\herm} \Mat{B}_{\fast}^{\pinv, \herm} \Mat{\Psi}\) using the type-II direct solver and NUFFT.}
    \STATE{Compress \(\Mat{B}_{\fast}^{\pinv} \Mat{A}\) into an HSS matrix \(\Mat{H}_{\HSS}\) using \(\{\Mat{\Omega}, \Mat{Y}, \Mat{\Psi}, \Mat{Z}\}\).}
    \STATE{Compute the URV factorization of \(\Mat{H}_{\HSS}\).}
    \STATE{Calculate \(\Mat{c} = \Mat{H}_{\HSS}^{-1} \Mat{B}_{\fast}^{\pinv} \Mat{f}\).}
  \end{algorithmic}
\end{algorithm}

Suppose the type-III NUDFT matrix \(\Mat{A} \in \complex^{M \times N}\) is defined by \(\Mat{A}(j, k) = \e^{-2 \pi \imag x_{j} \omega_{k}}\) for \(x_{j} \in [0, 1)\) and \(\omega_{k} \in [-1/2, N - 1 / 2)\).
Let \(\Mat{B}\) be the type-II NUDFT matrix defined by \(\Mat{B}(j, \ell) = \e^{-2 \pi \imag x_{j} \ell}\) where \(\ell \in \indexset{I}_{N; \sigma}\) is the index set defined in section~\ref{subsec:inudft3_motivations}.
Since \(\Mat{B}\) is a Vandermonde matrix with distinct nodes, it has full column rank and \(\Mat{B}^{\pinv} \Mat{B} = \Mat{I}\).
Multiplying \(\Mat{B}^{\pinv}\) on the left of~\eqref{eq:nudft3_matrix_approx} and denoting \(\Mat{H} = \Mat{B}^{\pinv} \Mat{A}\) gives
\begin{equation} \label{eq:nudft3_H}
    \Mat{H} = \Mat{Q} + \Mat{B}^{\pinv} \Mat{E}.
\end{equation}
Moreover, by multiplying \(\Mat{B}\) on the left and rearranging the equation, we have
\begin{equation} \label{eq:nudft3_BH}
    \Mat{A} = \Mat{B} \Mat{H} + \bigl(\Mat{I} - \Mat{B}\Mat{B}^{\pinv}\bigr) \Mat{E}.
\end{equation}

The two equations~\eqref{eq:nudft3_H} and~\eqref{eq:nudft3_BH} serve as the foundation of our algorithm.
From the discussion in section~\ref{subsec:inudft3_motivations}, \(\Mat{Q}\) can be well approximated by an HSS matrix.
Consequently, the matrix \(\Mat{H}\) in~\eqref{eq:nudft3_H}, viewed as a perturbation of \(\Mat{Q}\), can also be empirically compressed into an HSS matrix.
On the other hand, from~\eqref{eq:nudft3_BH}, \(\Mat{B} \Mat{H}\) can be regarded as an approximation of \(\Mat{A}\).
Noting that \(\Mat{I} - \Mat{B}\Mat{B}^{\pinv}\) is a projection matrix, the following results are direct consequences based on the discussion in section~\ref{subsec:inudft3_motivations}:
\begin{equation} \label{eq:nudft3_fast_approx_error}
    \begin{gathered}
        \expectation \|\bigl(\Mat{I} - \Mat{B}\Mat{B}^{\pinv}\bigr) \Mat{E}\|_{\fro}^2 \leq \expectation \|\Mat{E}\|_{\fro}^2 \leq \frac{2 \sin^2 (\pi \alpha)}{\pi^2 ((\sigma - 1) N / 2 + 1 / 2 - \alpha)} \|\Mat{A}\|_{\fro}^2.
    \end{gathered}
\end{equation}

Therefore, we approximate the NUDFT matrix \(\Mat{A}\) by the product of \(\Mat{B}\) and \(\Mat{H}\), where \(\Mat{B}\) is a type-II NUDFT matrix and \(\Mat{H}\) can be approximated by an HSS matrix.
In the preliminary stage, we construct a direct solver \(\Mat{B}_{\fast}\) for the type-II NUDFT matrix \(\Mat{B}\) as in Algorithm~\ref{alg:inudft2}.
Thus, the main task here is to find an efficient way to compress \(\Mat{H}\) into its HSS form.
This can be accomplished based on the following two observations:
\begin{itemize}
    \item The matrix-vector multiplications \(\Mat{A} \Mat{c}\) and \(\Mat{A}^{\herm} \Mat{f}\) correspond to the forward NUDFT computation, which can be performed very efficiently.
    \item The pseudoinverse (as well as the conjugate transpose of the pseudoinverse) of \(\Mat{B}_{\fast}\) can be applied rapidly by the fast solver introduced in section~\ref{subsec:inudft2}.
\end{itemize}
Therefore, by applying \(\Mat{B}_{\fast}^{\pinv} \Mat{A}\) on random test matrices, the black-box algorithm in section~\ref{subsec:hbs_black_box_construction} is used to obtain an HSS approximation of \(\Mat{H}\).
After construction, the fast structure in the type-III NUDFT problem can be expressed as
\begin{equation} \label{eq:nudft3_fast_forward}
    \Mat{A}_{\fast} = \Mat{B}_{\fast} \Mat{H}_{\HSS}.
\end{equation}

Once the HSS matrix \(\Mat{H}_{\HSS}\) has been prepared, we use the URV factorization to give an inversion of it.
Consequently, the solution of the least-squares problem is quite straightforward.
We directly compute \(\Mat{c} = \Mat{H}_{\HSS}^{-1} \Mat{B}_{\fast}^{\pinv} \Mat{f}\), which involves the solution processes for the type-II NUDFT matrix \(\Mat{B}\) and the HSS matrix \(\Mat{H}\).

This method can also serve as a preconditioner for iterative methods when the matrix \(\Mat{A}\) is too ill-conditioned.
When \(\Mat{A}\) is not too ill-conditioned, one can obtain a solution by applying PCG to the following normal equation:
\begin{equation} \label{eq:nudft_normal_eq}
    \Mat{A}^{\herm} \Mat{A} \Mat{c} = \Mat{A}^{\herm} \Mat{f}.
\end{equation}
Because of the fast matrix-vector multiplication, each iteration of CG can be computed efficiently.
In this circumstance, one can use \(\Mat{A}_{\fast}^{\pinv} \Mat{A}_{\fast}^{\pinv \herm}\) as a preconditioner of system~\eqref{eq:nudft_normal_eq}.

\subsection{Complexities} \label{subsec:complexity}

We analyze the complexity of Algorithm~\ref{alg:inudft3} under the assumption that the HSS rank of both \(\Mat{H}_{\HSS}\) and \(\Mat{\tilde{B}}_{\HSS}\) is \(\bigO(\log N)\).
The construction stage includes three steps:
As shown in section~\ref{subsec:inudft2}, the preliminary stage of constructing the type-II direct solver has a complexity of \(\bigO(M \log^{2} N)\).
The construction stage applies \(\Mat{A}\) and \(\Mat{B}_{\fast}^{\pinv}\) on random test matrices with \(\bigO(\log N)\) columns and recovers \(\Mat{H}_{\HSS}\) using the test matrices and sample matrices, which costs \(\bigO(M \log^{2} N  + N \log^{2} N)\) and \(\bigO(N \log^{2} N)\).
Therefore, the total complexity of the construction stage is \(\bigO(M \log^{2} N)\).
The factorization stage involves computing the URV factorization of \(\Mat{H}_{\HSS}\), which has a cost of \(\bigO(N \log^{2} N)\).
The solution stage includes the applications of \(\Mat{B}_{\fast}^{\pinv}\) and \(\Mat{H}_{\HSS}^{\pinv}\), with complexities of \(\bigO(M \log N)\) and \(\bigO(N \log N)\), respectively.
Hence, the total complexity for solution is \(\bigO(M \log N)\).

\section{Numerical results} \label{sec:numerical_results}

In this section, we conduct several tests on our proposed algorithm.
Throughout all experiments, two cases of type-III NUDFT problem are considered:
\begin{itemize}
    \item Perturbation-perturbation:
    Sample points and frequencies are perturbed from uniform points and continuous integers, respectively, i.e., \(x_{j} = (j + \beta \phi_{j}) / M\) and \(\omega_{k} = k + \alpha \psi_{k}\) where each \(\phi_{j}\) and \(\psi_{k}\) is drawn from a uniform distribution on \([-1, 1]\) and \(0 \leq \alpha, \beta < 1 / 2\) control the nonuniformity.
    \item Random-perturbation: Sample points are i.i.d.\ uniform random variables on \([0, 1)\), and the frequencies are perturbed integers given by \(\omega_{k} = k + \alpha \psi_{k}\).
\end{itemize}
We briefly explain our experimental settings.
The number of frequency points \(N\) is chosen to be \(2^{p}\) for \(p = 16 : 20\), and the number of sample points \(M\) is set to be \(4 N\).
The leaf size of each HSS matrix is fixed at \(256\), the accuracy parameter in the fADI construction of \(\Mat{\tilde{B}}_{\HSS}\) is set to be \(10^{-6}\), and the HSS rank of \(\Mat{H}_{\HSS}\) in~\eqref{eq:nudft3_fast_forward} is set to be approximately \(10 \log_{2} (N)\) for its black-box construction.
The exact coefficient vector \(\Mat{c}_{\exact}\) is generated by a random vector in \(\complex^{N}\), and the target values are obtained by applying the type-III NUFFT to \(\Mat{c}_{\exact}\).
For the perturbation-perturbation case, \(\beta\) is fixed to be \(0.4\) and \(\alpha\) is chosen to be \(0.1\) or \(0.4\).
For the random-perturbation case, \(\alpha\) is set to be \(0.1\) and \(0.4\).
For iterative methods, we use the MATLAB \texttt{lsqr} command~\cite{Paige_Saunders_1982}, with tolerance set to be \(10^{-12}\) and the maximum number of iterations set to be \(500\).
We denote the iterative method without and with preconditioner as ``lsqr'' and ``plsqr,'' respectively.
All algorithms are implemented in MATLAB R2023b where the MATLAB FFT library and the FINUFFT library~\cite{Barnett_Magland_Klinteberg_2019} are utilized for all corresponding fast matrix-vector multiplications.
All experiments are carried out on a server with Intel Gold 6226R CPU at 2.90 GHz and 1000.6 GB of RAM.

The following values are taken into consideration when evaluating the algorithm:
We use \(t_{\cB}\), \(t_{\fB}\), \(t_{\cH}\), and \(t_{\fH}\) to denote the time cost for constructing \(\Mat{\tilde{B}}_{\HSS}\), factorizing \(\Mat{\tilde{B}}_{\HSS}\), constructing \(\Mat{H}_{\HSS}\), and factorizing \(\Mat{H}_{\HSS}\), respectively.
For each solution, its relative residual and relative error are defined as \(r_{\solve} = \|\Mat{A} \Mat{c} - \Mat{f}\|_{2} / \|\Mat{f}\|_{2}\) and \(e_{\solve} = \|\Mat{c} - \Mat{c}_{\exact}\|_{2} / \|\Mat{c}_{\exact}\|_{2}\), respectively.
For direct methods, we use \(t_{\direct}\) to denote the time cost for solving the system, while for iterative methods, we use \(t_{\pre}\) and \(t_{\iter}\) to denote the time cost in the preparation and iteration stage, respectively.
More specifically, \(t_{\pre} = t_{\cB} + t_{\fB} + t_{\cH} + t_{\fH}\), and \(t_{\iter}\) only includes the iteration time of lsqr or plsqr.
Additionally, \(n_{\iter}\) is the number of iterations to solve~\eqref{eq:nudft_normal_eq} using lsqr or plsqr to a tolerance of \(10^{-12}\), with a maximum number of iterations set to be \(500\).

\subsection{Exploring properties of the type-III NUDFT} \label{subsec:properties_nudft3}

\begin{figure}[htbp]
    \centering
    \begin{subfigure}{0.35\textwidth}
        \centering
        \includegraphics[width=\linewidth]{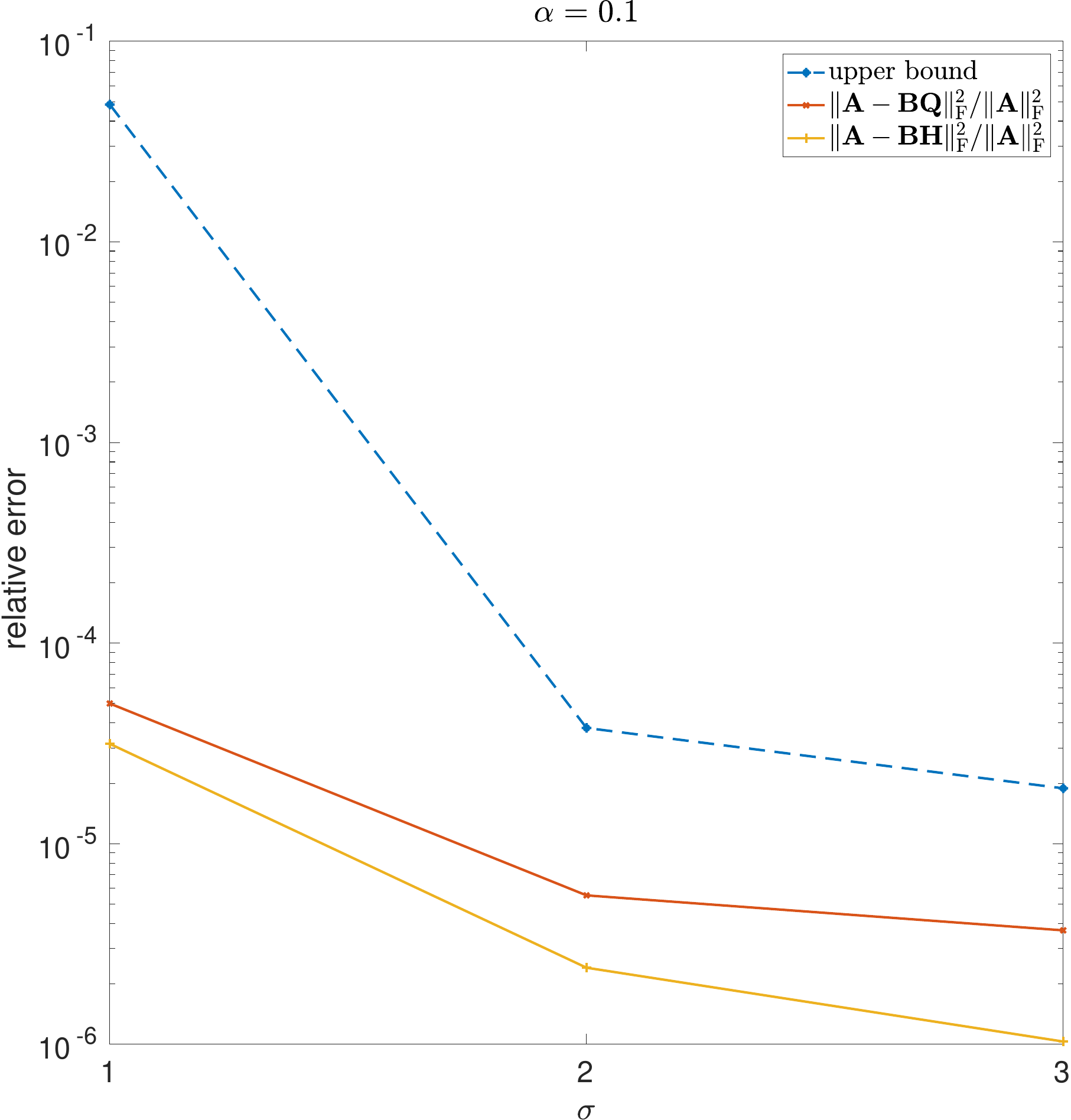}
    \end{subfigure}
    \hspace{0.01\textwidth}
    \begin{subfigure}{0.35\textwidth}
        \centering
        \includegraphics[width=\linewidth]{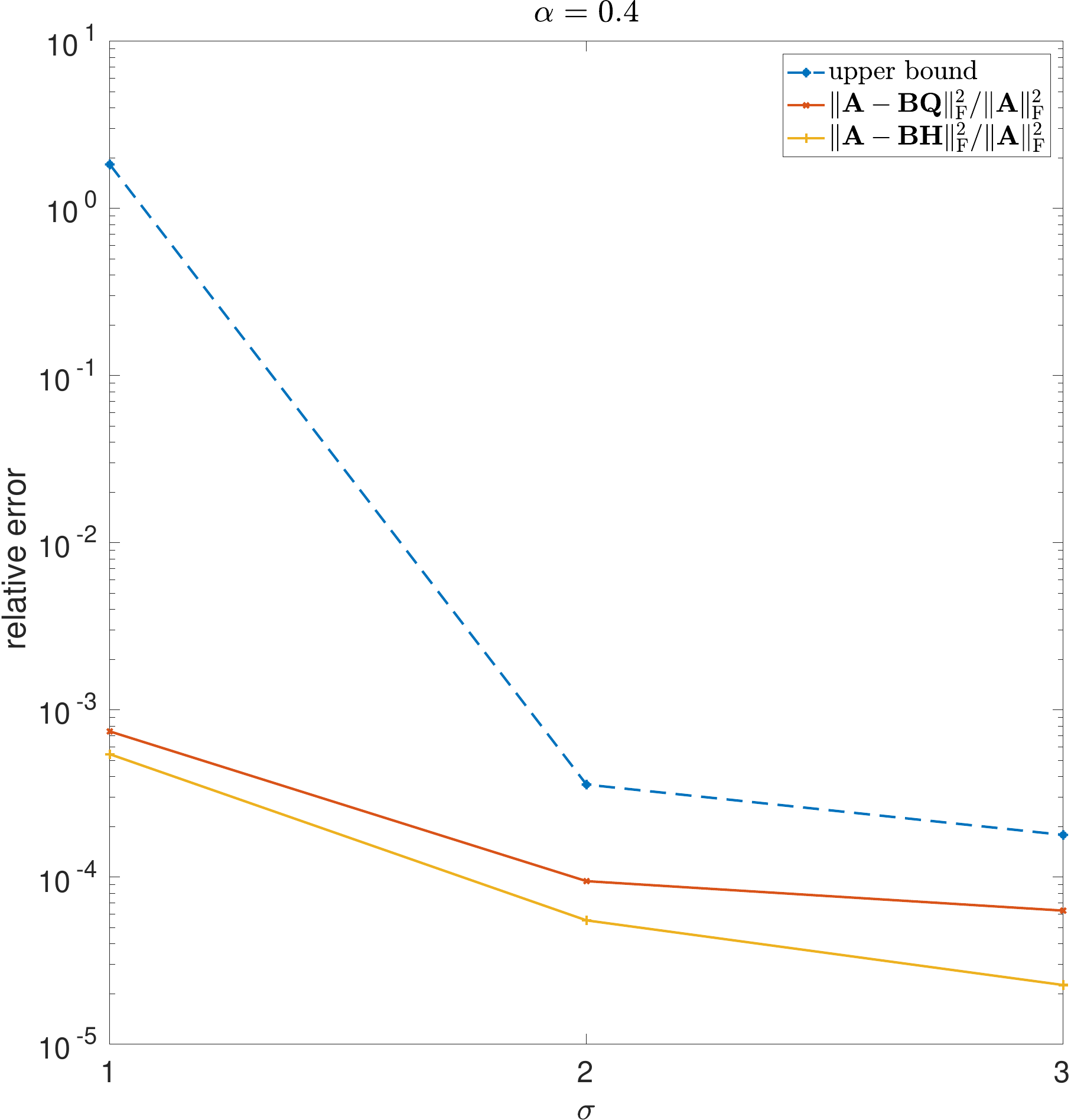}
    \end{subfigure}
    \caption{Approximation errors for various \(\sigma\) and \(\alpha\).
    Left: \(\alpha = 0.1\).
    Right: \(\alpha = 0.4\).
    The upper bound is given by Proposition~\ref{prop:nudft3_matrix_approx_error} and~\eqref{eq:nudft3_fast_approx_error}.}
    \label{fig:typeIII_approx}
\end{figure}

We verify the approximation error in Proposition~\ref{prop:nudft3_matrix_approx_error} and~\eqref{eq:nudft3_fast_approx_error} by approximating a \(4096 \times 1024\) type-III NUDFT matrix \(\Mat{A}\) under the random-perturbation case.
The approximation error is defined as \(\|\Mat{A} - \Mat{B} \Mat{Q}\|_{\fro}^{2} / \|\Mat{A}\|_{\fro}^{2}\) and \(\|\Mat{A} - \Mat{B} \Mat{H}\|_{\fro}^{2} / \|\Mat{A}\|_{\fro}^{2}\), respectively.
Both errors are estimated by an average over \(5\) simulations.
Figure~\ref{fig:typeIII_approx} shows the error for different values of \(\sigma\) and \(\alpha\).
It can be observed that the approximation error decreases as \(\sigma\) increases and as \(\alpha\) decreases, which is consistent with our intuitions and the theoretical results.

\subsection{Perturbation-perturbation case} \label{subsec:perturb_perturb}

In~\cite{Austin_Trefethen_2017}, the authors discuss the convergence of interpolation and quadrature for periodic functions when using perturbed points.
In~\cite{Yu_Townsend_2023}, the authors demonstrate that the condition number of a tall-and-skinny type-II NUDFT matrix with perturbed sample points can be bounded by a constant that depends only on the perturbation size \(\beta\) and the ratio \(M / N\).
Although there are no directly related generalizations for the type-III case, this finding provides useful insights.
When the perturbation \(\alpha\) of frequencies is small, we expect similar behaviors of the NUDFT matrix.

\begin{figure}[tbhp]
    \centering
    \begin{subfigure}{0.35\textwidth}
        \centering
        \includegraphics[width=\linewidth]{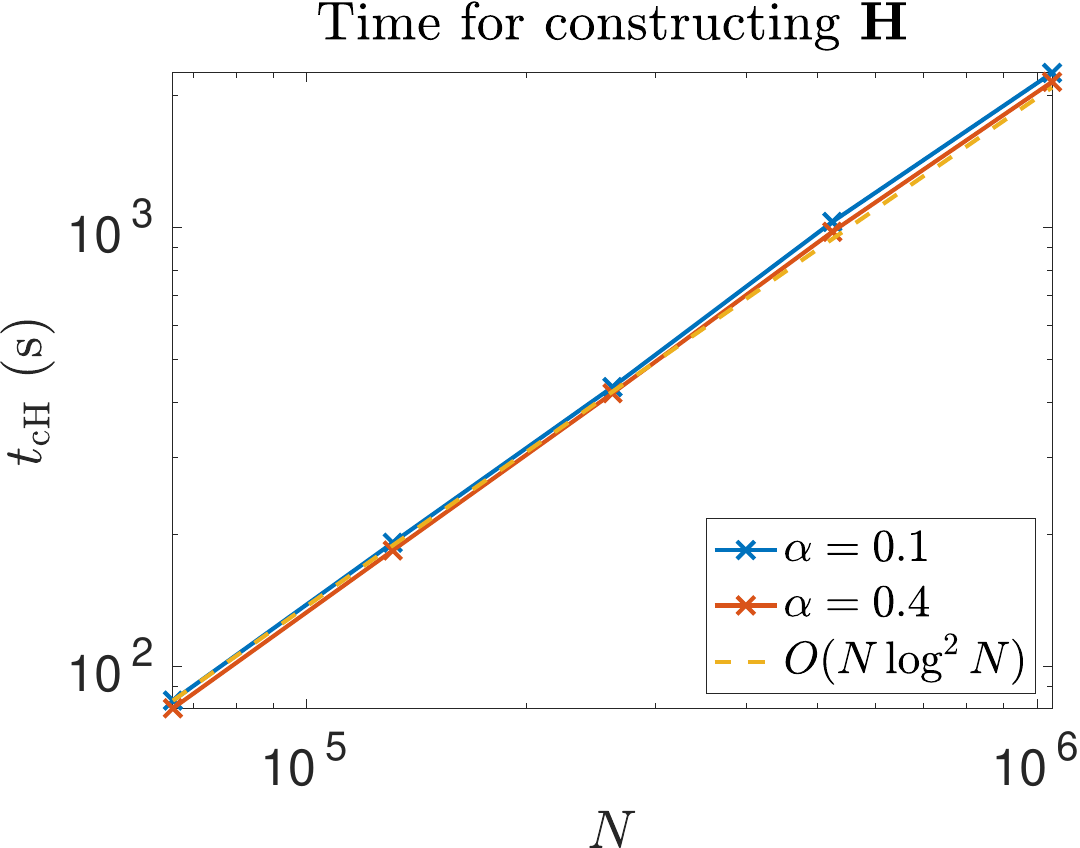}
        \caption{}
        \label{fig:pp_t_cH}
    \end{subfigure}
    \hspace{0.01\textwidth}
    \begin{subfigure}{0.35\textwidth}
        \centering
        \includegraphics[width=\linewidth]{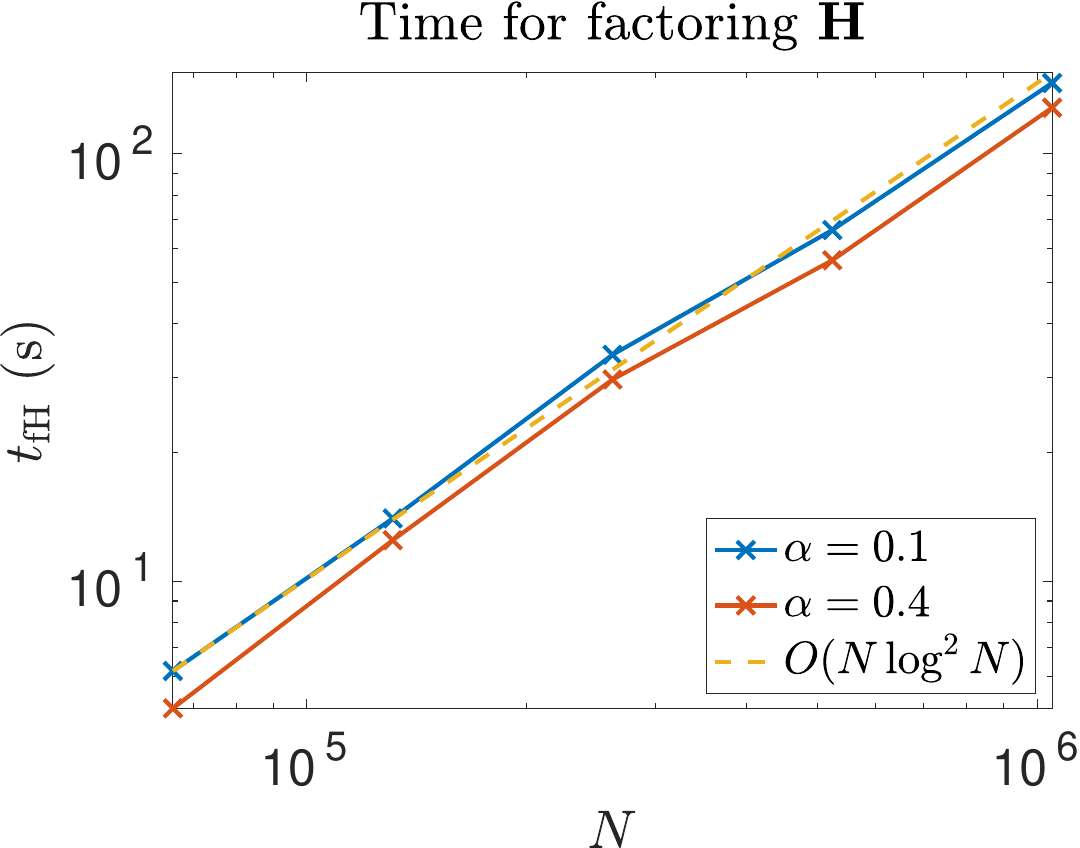}
        \caption{}
        \label{fig:pp_t_fH}
    \end{subfigure}
    \begin{subfigure}{0.35\textwidth}
        \centering
        \includegraphics[width=\linewidth]{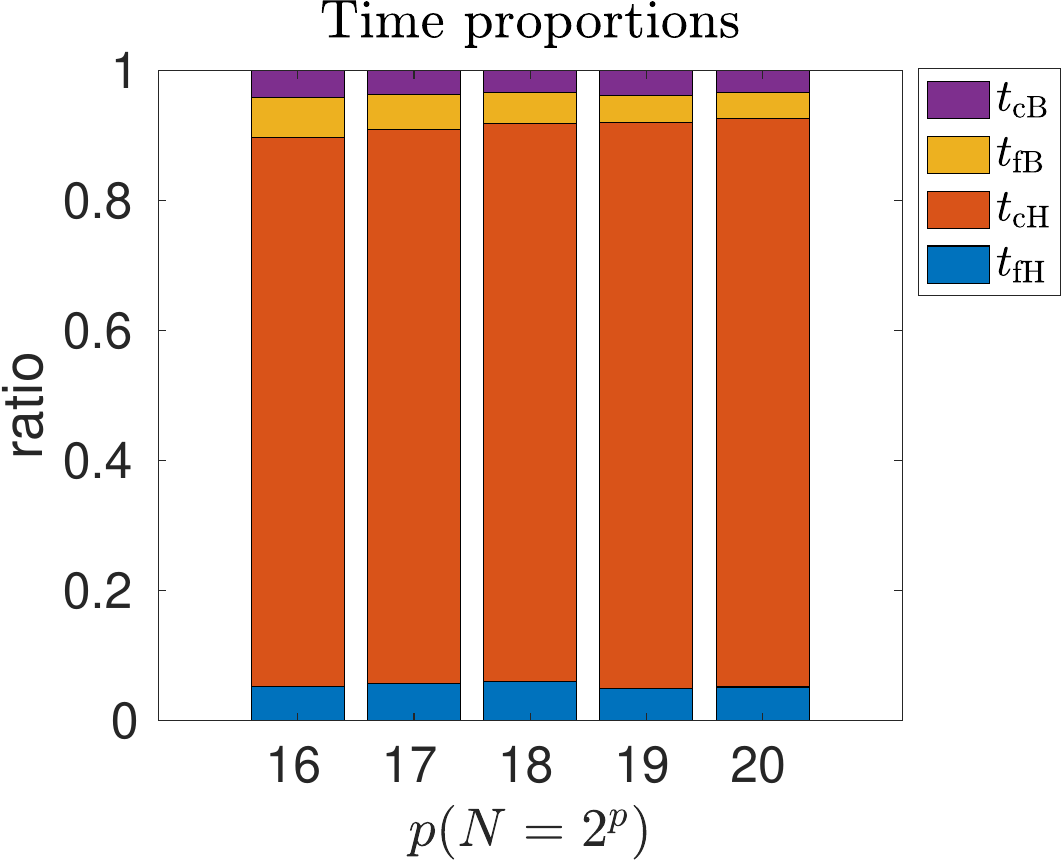}
       \caption{}
        \label{fig:pp_ratio}
    \end{subfigure}
    \begin{subfigure}{0.35\textwidth}
        \centering
        \includegraphics[width=\linewidth]{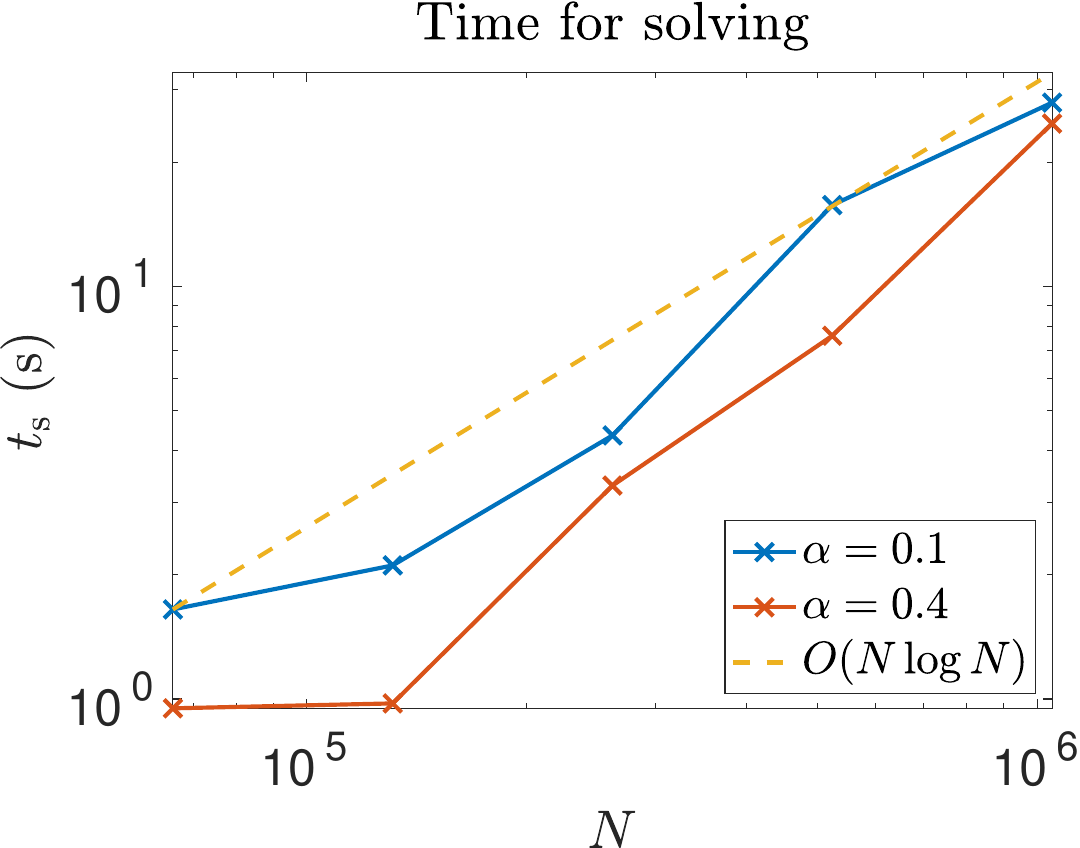}
        \caption{}
        \label{fig:pp_t_s}
    \end{subfigure}
    \begin{subfigure}{0.35\textwidth}
        \centering
        \includegraphics[width=\linewidth]{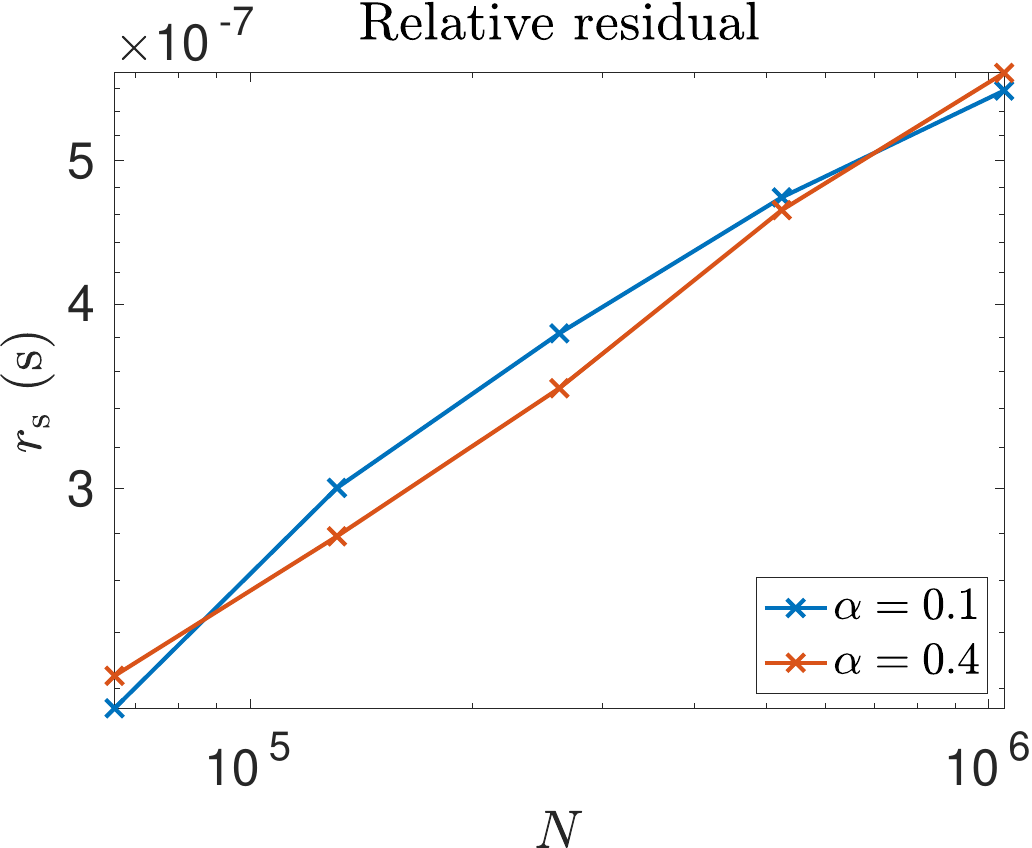}
       \caption{}
        \label{fig:pp_r_s}
    \end{subfigure}
    \begin{subfigure}{0.35\textwidth}
        \centering
        \includegraphics[width=\linewidth]{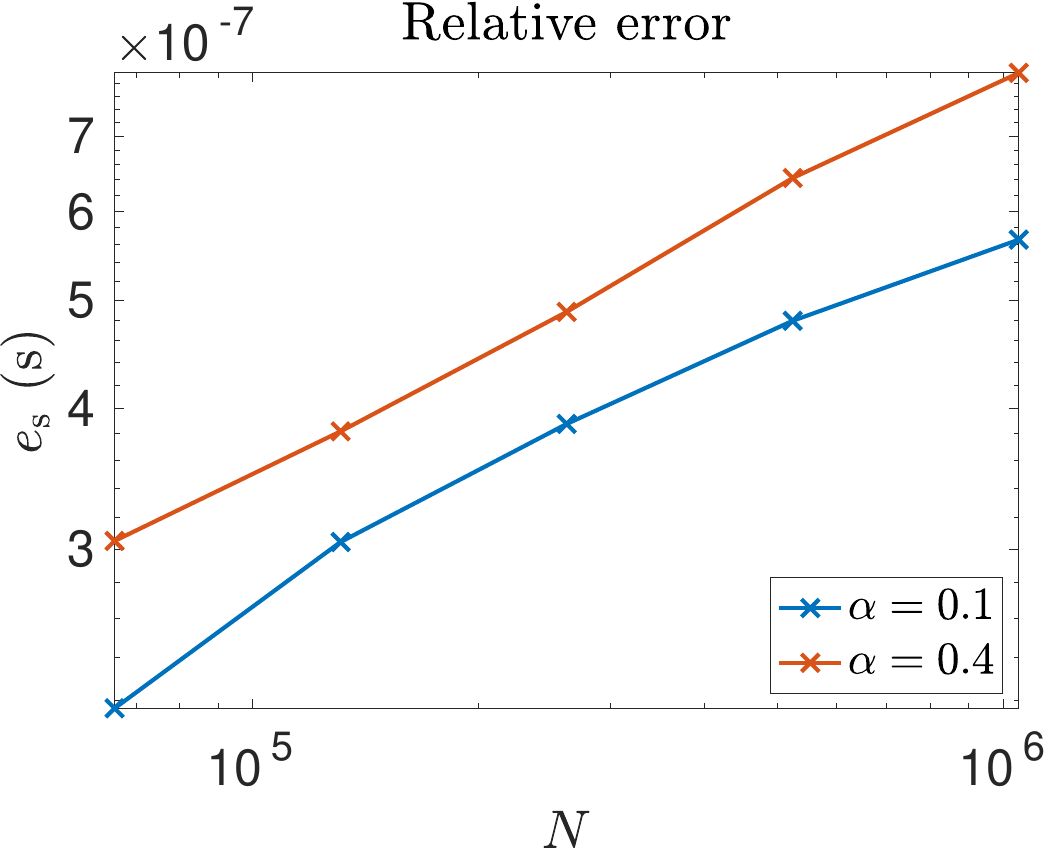}
       \caption{}
        \label{fig:pp_e_s}
    \end{subfigure}
    \caption{Results of the direct solver under the perturbation-perturbation case for \(\alpha = 0.1\) and \(0.4\) and \(\beta = 0.4\).}
    \label{fig:pp_direct}
\end{figure}

\begin{table}[tbhp]
\centering
\begin{tabular}{c cccc ccc}
\toprule
\(N\) & \(t_{\cB}\) & \(t_{\fB}\) & \(t_{\cH}\) & \(t_{\fH}\) & \(t_{\direct}\) & \(r_{\s}\) & \(e_{\s}\) \\
\midrule
\(2^{16}\) & 4.1e+00 & 6.0e+00 & 8.4e+01 & 6.2e+00 & 1.6e+00 & 2.1e-07 & 2.2e-07 \\
\midrule
\(2^{17}\) & 8.1e+00 & 1.1e+01 & 1.9e+02 & 1.4e+01 & 2.1e+00 & 3.0e-07 & 3.0e-07 \\
\midrule
\(2^{18}\) & 1.6e+01 & 2.3e+01 & 4.3e+02 & 3.4e+01 & 4.4e+00 & 3.8e-07 & 3.9e-07 \\
\midrule
\(2^{19}\) & 4.0e+01 & 4.6e+01 & 1.0e+03 & 6.6e+01 & 1.6e+01 & 4.7e-07 & 4.8e-07 \\
\midrule
\(2^{20}\) & 9.5e+01 & 9.7e+01 & 2.3e+03 & 1.5e+02 & 2.8e+01 & 5.6e-07 & 5.7e-07 \\
\bottomrule
\end{tabular}
\caption{Time and error results of the direct solver under the perturbation-perturbation case, \(\alpha = 0.1\) and \(\beta = 0.4\).}
\label{tab:direct_pp_1e-1_4e-1}
\end{table}

\begin{table}[tbhp]
\centering
\begin{tabular}{c cccc ccc}
\toprule
\(N\) & \(t_{\cB}\) & \(t_{\fB}\) & \(t_{\cH}\) & \(t_{\fH}\) & \(t_{\direct}\) & \(r_{\s}\) & \(e_{\s}\) \\
\midrule
\(2^{16}\) & 3.9e+00 & 5.8e+00 & 8.0e+01 & 5.1e+00 & 9.5e-01 & 2.2e-07 & 3.1e-07 \\
\midrule
\(2^{17}\) & 7.9e+00 & 1.1e+01 & 1.8e+02 & 1.2e+01 & 9.8e-01 & 2.8e-07 & 3.8e-07 \\
\midrule
\(2^{18}\) & 1.6e+01 & 2.3e+01 & 4.2e+02 & 3.0e+01 & 3.3e+00 & 3.5e-07 & 4.9e-07 \\
\midrule
\(2^{19}\) & 4.3e+01 & 4.7e+01 & 9.8e+02 & 5.6e+01 & 7.6e+00 & 4.6e-07 & 6.4e-07 \\
\midrule
\(2^{20}\) & 8.0e+01 & 9.8e+01 & 2.2e+03 & 1.3e+02 & 2.5e+01 & 5.7e-07 & 8.0e-07 \\
\bottomrule
\end{tabular}
\caption{Time and error results of the direct solver under the perturbation-perturbation case, \(\alpha = 0.4\) and \(\beta = 0.4\).}
\label{tab:direct_pp_4e-1_4e-1}
\end{table}

\begin{table}[tbhp]
\centering
\begin{tabular}{c cccc ccc}
\toprule
\(N\) & Method & \(t_{\pre}\) & \(t_{\iter}\) & \(n_{\iter}\) & \(r_{\s}\) & \(e_{\s}\) \\
\midrule
\multirow{2}{*}{\(2^{16}\)} & lsqr & - & 7.8e-01 & 17 & 7.4e-13 & 8.0e-13 \\
 & plsqr & 1.0e+02 & 6.9e+00 & 5 & 4.1e-13 & 4.3e-13 \\
\midrule
\multirow{2}{*}{\(2^{17}\)} & lsqr & - & 1.5e+00 & 17 & 6.7e-13 & 7.1e-13 \\
 & plsqr & 2.2e+02 & 1.2e+01 & 5 & 4.7e-14 & 4.9e-14 \\
\midrule
\multirow{2}{*}{\(2^{18}\)} & lsqr & - & 3.1e+00 & 17 & 6.4e-13 & 6.9e-13 \\
 & plsqr & 5.1e+02 & 2.6e+01 & 5 & 2.8e-13 & 3.0e-13 \\
\midrule
\multirow{2}{*}{\(2^{19}\)} & lsqr & - & 7.4e+00 & 17 & 9.3e-13 & 1.0e-12 \\
 & plsqr & 1.2e+03 & 6.1e+01 & 5 & 1.3e-13 & 1.4e-13 \\
\midrule
\multirow{2}{*}{\(2^{20}\)} & lsqr & - & 1.9e+01 & 17 & 8.0e-13 & 8.7e-13 \\
 & plsqr & 2.6e+03 & 1.3e+02 & 5 & 5.6e-14 & 5.7e-14 \\
\bottomrule
\end{tabular}
\caption{Time and error results of the iterative solver under the perturbation-perturbation case, \(\alpha = 0.1\) and \(\beta = 0.4\).}
\label{tab:iterative_pp_1e-1_4e-1}
\end{table}

\begin{table}[tbhp]
\centering
\begin{tabular}{c cccc ccc}
\toprule
\(N\) & Method & \(t_{\pre}\) & \(t_{\iter}\) & \(n_{\iter}\) & \(r_{\s}\) & \(e_{\s}\) \\
\midrule
\multirow{2}{*}{\(2^{16}\)} & lsqr & - & 4.1e+00 & 95 & 8.8e-13 & 2.9e-12 \\
 & plsqr & 9.5e+01 & 8.5e+00 & 8 & 4.0e-13 & 7.3e-13 \\
\midrule
\multirow{2}{*}{\(2^{17}\)} & lsqr & - & 8.6e+00 & 99 & 1.0e-12 & 4.0e-12 \\
 & plsqr & 2.2e+02 & 2.0e+01 & 9 & 3.3e-15 & 6.1e-15 \\
\midrule
\multirow{2}{*}{\(2^{18}\)} & lsqr & - & 1.7e+01 & 101 & 7.8e-13 & 2.9e-12 \\
 & plsqr & 4.9e+02 & 4.7e+01 & 9 & 6.6e-14 & 1.2e-13 \\
\midrule
\multirow{2}{*}{\(2^{19}\)} & lsqr & - & 3.4e+01 & 104 & 8.9e-13 & 3.5e-12 \\
 & plsqr & 1.1e+03 & 9.7e+01 & 9 & 1.6e-14 & 3.2e-14 \\
\midrule
\multirow{2}{*}{\(2^{20}\)} & lsqr & - & 1.2e+02 & 105 & 9.9e-13 & 4.0e-12 \\
 & plsqr & 2.5e+03 & 2.1e+02 & 9 & 3.6e-14 & 6.5e-14 \\
\bottomrule
\end{tabular}
\caption{Time and error results of the iterative solver under the perturbation-perturbation case, \(\alpha = 0.4\) and \(\beta = 0.4\).}
\label{tab:iterative_pp_4e-1_4e-1}
\end{table}

Figure~\ref{fig:pp_direct} and Tables~\ref{tab:direct_pp_1e-1_4e-1} and~\ref{tab:direct_pp_4e-1_4e-1} show the results of the direct solver under the perturbation-perturbation case for \(\alpha = 0.1\) and \(0.4\).
It can be observed that the time costs of different components in the algorithm scale as predicted in section~\ref{subsec:complexity}, except for the solution stage.
This is because the time cost of the solution stage is quite small compared to others.
As we expect, the runtimes for \(\alpha = 0.4\) and \(\alpha = 0.1\) are quite similar, since our algorithm is not sensitive to distribution of sample points and frequencies.
In all examples, the most time-consuming part is the black-box construction of \(\Mat{H}_{\HSS}\), accounting for nearly \(90\%\) of the total time cost; this can also be inferred from the proportion of different components shown in Figure~\ref{fig:pp_direct}(\subref{fig:pp_ratio}).
To be more specific, when \(\alpha = 0.4\) and \(N = 2^{20}\), the runtimes of \(t_{\cB}\), \(t_{\fB}\), \(t_{\cH}\), \(t_{\fH}\) are \(80\)s, \(98\)s, \(2200\)s, and \(130\)s, respectively, while the solution only takes \(25\)s.
On the other hand, the relative residual and relative error are at the order of \(10^{-7}\) for all examples, meaning that our algorithm is quite accurate and stable.

The results of the iterative solver are shown in Tables~\ref{tab:iterative_pp_1e-1_4e-1} and~\ref{tab:iterative_pp_4e-1_4e-1}.
For both values of \(\alpha\), the NUDFT matrix is quite well-conditioned, since the number of iterations for lsqr is quite small.
When \(\alpha = 0.1\), the number of iterations is \(17\) for all examples, while when \(\alpha = 0.4\), the number of iterations is around \(100\).
On the other hand, our algorithm serves as a very efficient preconditioner for the iterative method, since the number of iterations for plsqr is at most \(5\) and \(9\) for \(\alpha = 0.1\) and \(0.4\), respectively.
However, the time cost for plsqr is higher than that of lsqr, since each iteration of plsqr is more expensive and the problem is easy for applying lsqr directly.

\subsection{Random-perturbation case} \label{subsec:random_perturb}

For the random-perturbation case, we did not find any literature discussing the conditioning of the type-III NUDFT matrix.
For the type-II case when sample points are i.i.d.\ uniform random variables, the paper~\cite{Kammerer_Ullrich_Volkmer_2021} contains a study of the conditioning of the NUDFT matrix, showing that it is well-conditioned with high probability when \(M \gtrsim N \log N\).
This implies that the random-perturbation case is more challenging than the perturbation-perturbation case.

\begin{figure}[tbhp]
    \centering
    \begin{subfigure}{0.35\textwidth}
        \centering
        \includegraphics[width=\linewidth]{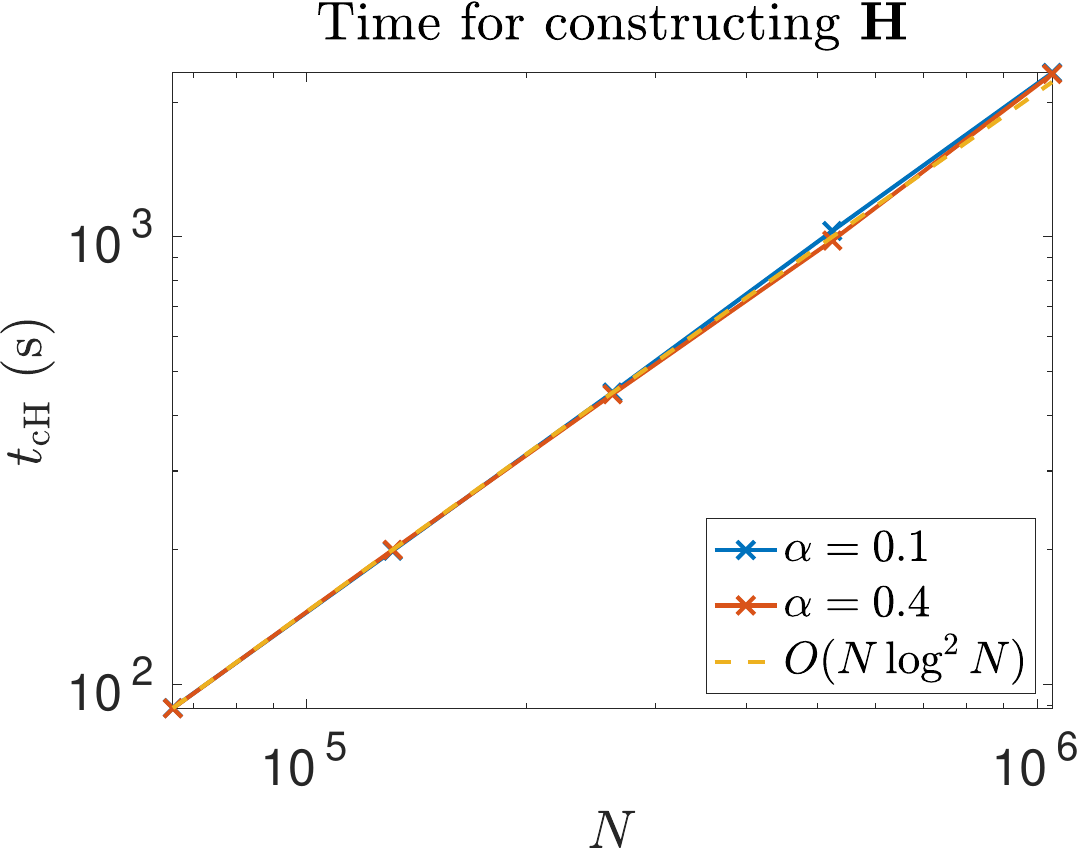}
        \caption{}
        \label{fig:rp_t_cH}
    \end{subfigure}
    \hspace{0.01\textwidth}
    \begin{subfigure}{0.35\textwidth}
        \centering
        \includegraphics[width=\linewidth]{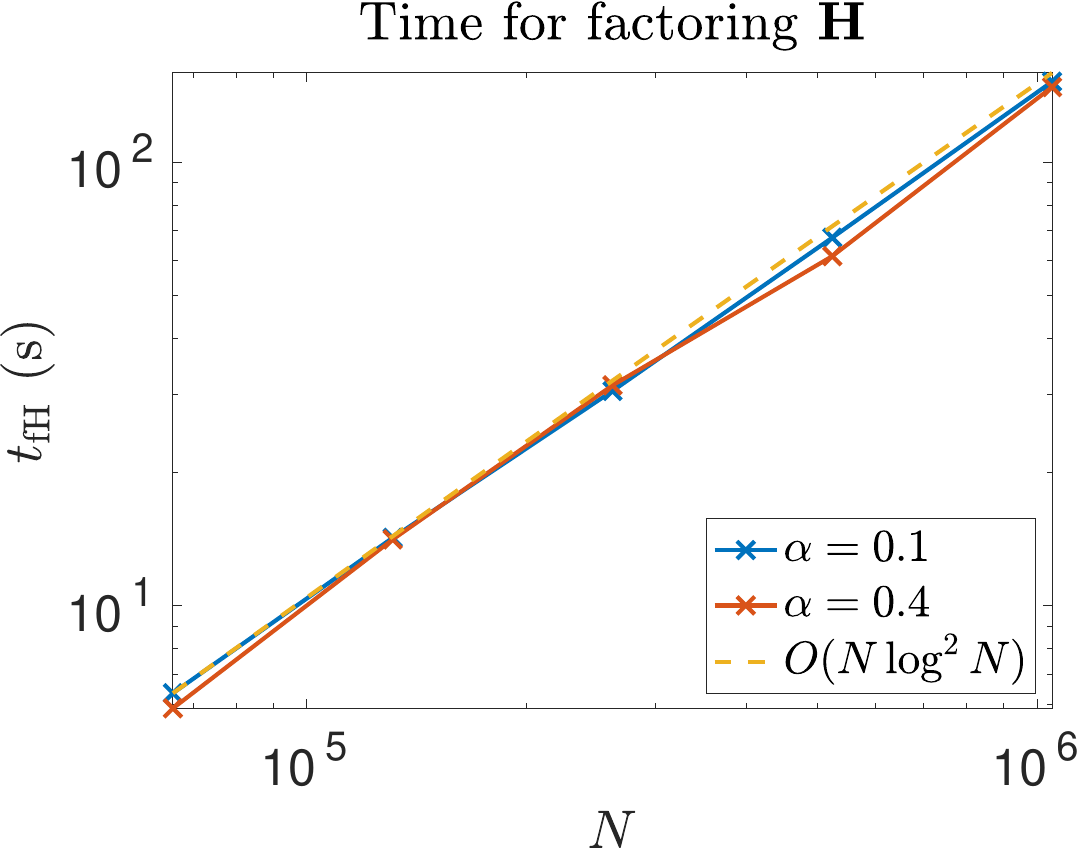}
        \caption{}
        \label{fig:rp_t_fH}
    \end{subfigure}
    \begin{subfigure}{0.35\textwidth}
        \centering
        \includegraphics[width=\linewidth]{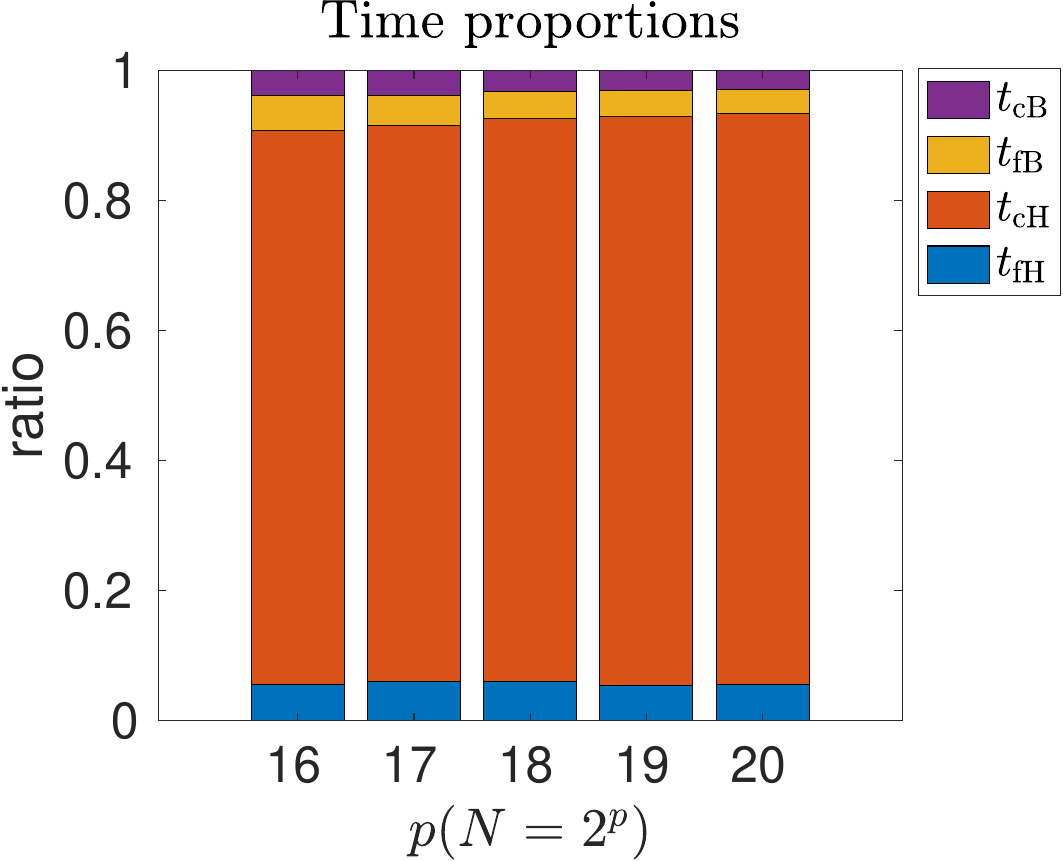}
        \caption{}
        \label{fig:rp_ratio}
    \end{subfigure}
    \begin{subfigure}{0.35\textwidth}
        \centering
        \includegraphics[width=\linewidth]{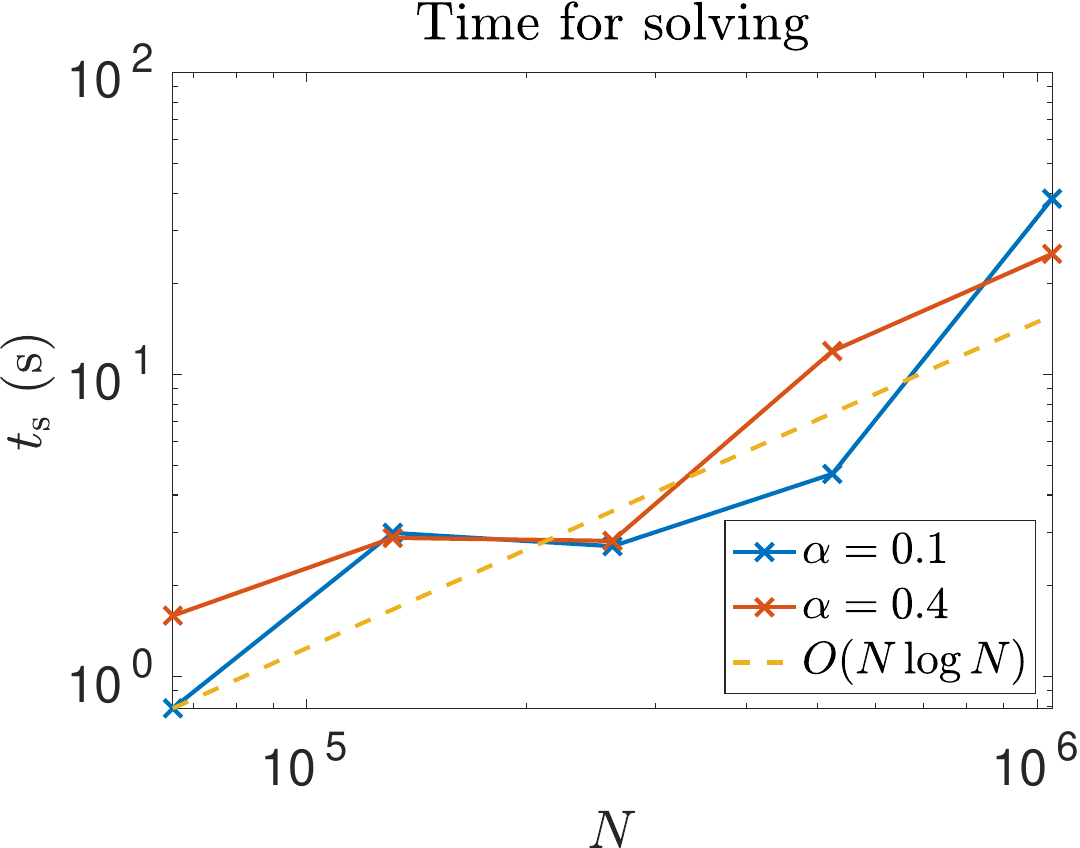}
        \caption{}
        \label{fig:rp_t_s}
    \end{subfigure}
    \begin{subfigure}{0.35\textwidth}
        \centering
        \includegraphics[width=\linewidth]{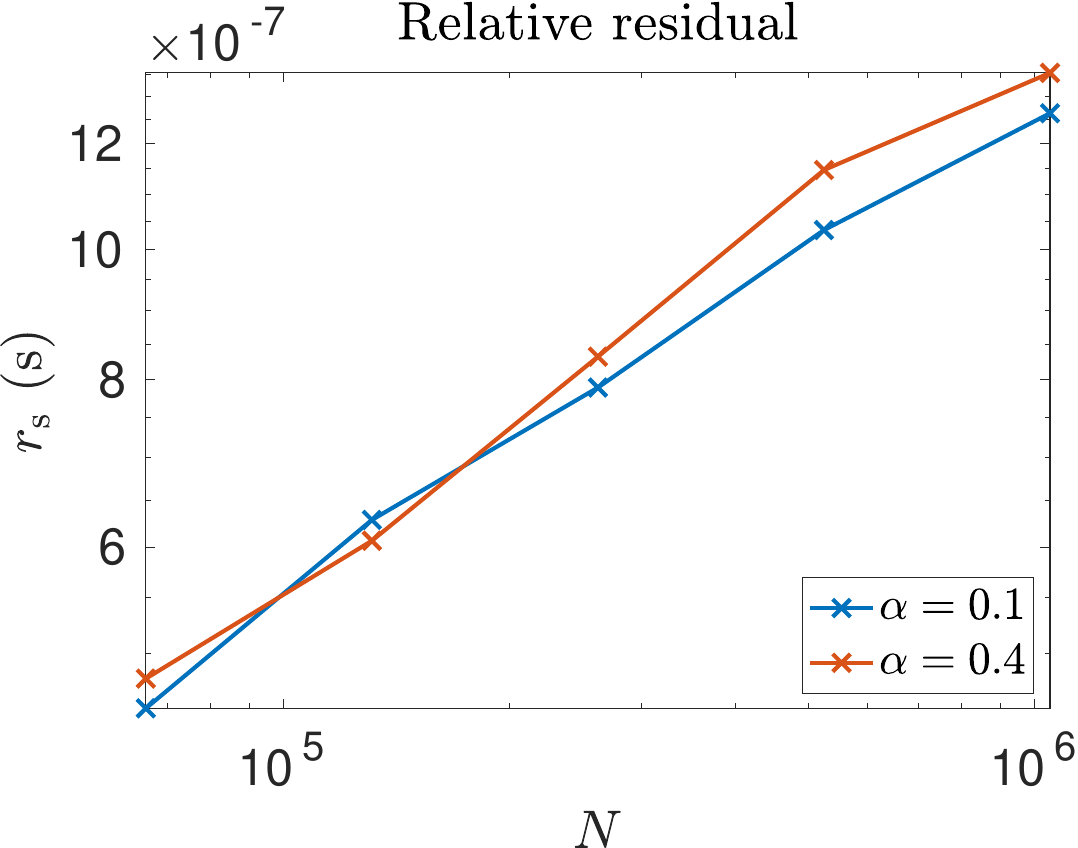}
        \caption{}
        \label{fig:rp_r_s}
    \end{subfigure}
    \begin{subfigure}{0.35\textwidth}
        \centering
        \includegraphics[width=\linewidth]{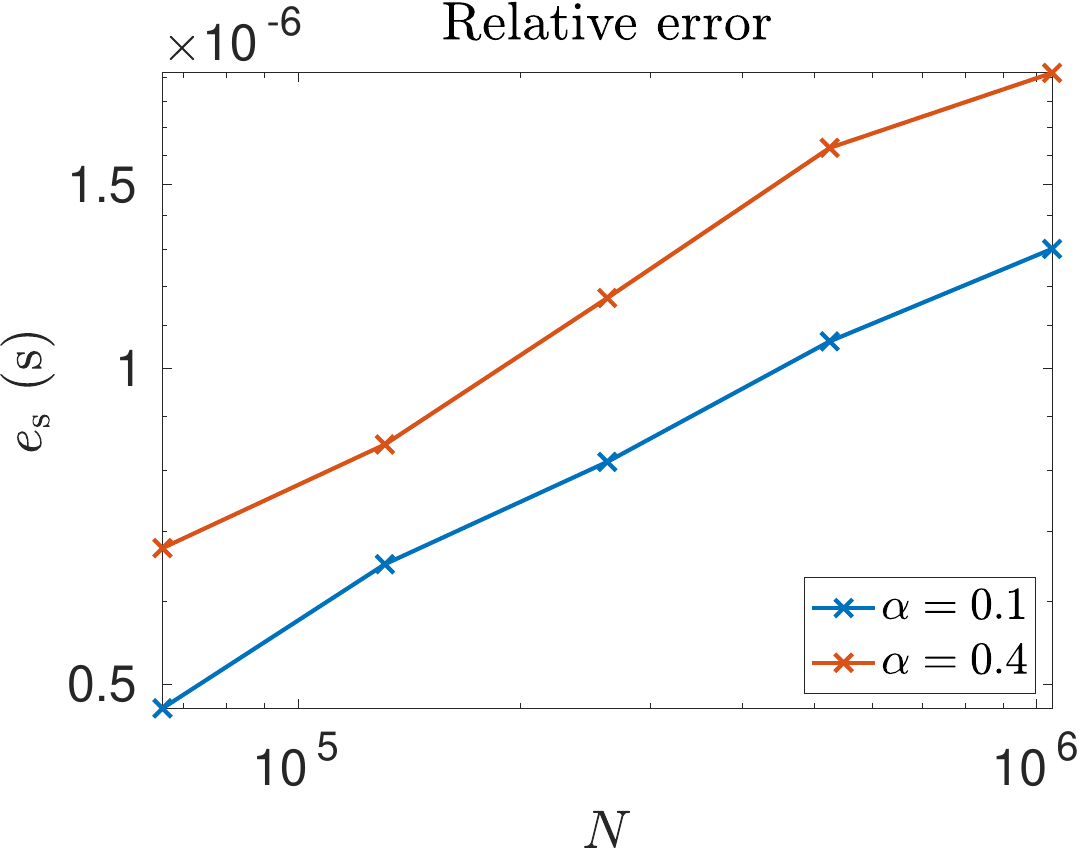}
        \caption{}
        \label{fig:rp_e_s}
    \end{subfigure}
    \caption{Results of the direct solver under the random-perturbation case for \(\alpha = 0.1\) and \(0.4\).}
    \label{fig:rp_direct}
\end{figure}

\begin{table}[tbhp]
\centering
\begin{tabular}{c cccc ccc}
\toprule
\(N\) & \(t_{\cB}\) & \(t_{\fB}\) & \(t_{\cH}\) & \(t_{\fH}\) & \(t_{\direct}\) & \(r_{\s}\) & \(e_{\s}\) \\
\midrule
\(2^{16}\) & 3.9e+00 & 5.4e+00 & 8.9e+01 & 6.4e+00 & 7.9e-01 & 4.5e-07 & 4.7e-07 \\
\midrule
\(2^{17}\) & 8.1e+00 & 1.1e+01 & 2.0e+02 & 1.4e+01 & 3.0e+00 & 6.3e-07 & 6.5e-07 \\
\midrule
\(2^{18}\) & 1.7e+01 & 2.2e+01 & 4.5e+02 & 3.1e+01 & 2.7e+00 & 7.9e-07 & 8.2e-07 \\
\midrule
\(2^{19}\) & 3.3e+01 & 4.4e+01 & 1.0e+03 & 6.8e+01 & 4.7e+00 & 1.0e-06 & 1.1e-06 \\
\midrule
\(2^{20}\) & 7.0e+01 & 9.1e+01 & 2.3e+03 & 1.5e+02 & 3.8e+01 & 1.3e-06 & 1.3e-06 \\
\bottomrule
\end{tabular}
\caption{Time and error results of the direct solver under the random-perturbation case, \(\alpha = 0.1\).}
\label{tab:direct_rp_1e-1}
\end{table}

\begin{table}[tbhp]
\centering
\begin{tabular}{c cccc ccc}
\toprule
\(N\) & \(t_{\cB}\) & \(t_{\fB}\) & \(t_{\cH}\) & \(t_{\fH}\) & \(t_{\direct}\) & \(r_{\s}\) & \(e_{\s}\) \\
\midrule
\(2^{16}\) & 3.9e+00 & 5.5e+00 & 8.9e+01 & 5.9e+00 & 1.6e+00 & 4.8e-07 & 6.7e-07 \\
\midrule
\(2^{17}\) & 9.0e+00 & 1.1e+01 & 2.0e+02 & 1.4e+01 & 2.9e+00 & 6.1e-07 & 8.5e-07 \\
\midrule
\(2^{18}\) & 1.6e+01 & 2.2e+01 & 4.5e+02 & 3.1e+01 & 2.8e+00 & 8.3e-07 & 1.2e-06 \\
\midrule
\(2^{19}\) & 3.4e+01 & 4.4e+01 & 9.8e+02 & 6.1e+01 & 1.2e+01 & 1.1e-06 & 1.6e-06 \\
\midrule
\(2^{20}\) & 7.7e+01 & 9.6e+01 & 2.3e+03 & 1.5e+02 & 2.5e+01 & 1.4e-06 & 1.9e-06 \\
\bottomrule
\end{tabular}
\caption{Time and error results of the direct solver under the random-perturbation case, \(\alpha = 0.4\).}
\label{tab:direct_rp_4e-1}
\end{table}

\begin{table}[tbhp]
\centering
\begin{tabular}{c cccc ccc}
\toprule
\(N\) & Method & \(t_{\pre}\) & \(t_{\iter}\) & \(n_{\iter}\) & \(r_{\s}\) & \(e_{\s}\) \\
\midrule
\multirow{2}{*}{\(2^{16}\)} & lsqr & - & 2.1e+01 & 500 & 2.3e-12 & 3.0e-11 \\
 & plsqr & 1.0e+02 & 9.2e+00 & 5 & 5.3e-14 & 2.0e-13 \\
\midrule
\multirow{2}{*}{\(2^{17}\)} & lsqr & - & 4.3e+01 & 500 & 5.9e-06 & 1.7e-03 \\
 & plsqr & 2.3e+02 & 1.4e+01 & 5 & 2.8e-13 & 1.6e-11 \\
\midrule
\multirow{2}{*}{\(2^{18}\)} & lsqr & - & 8.4e+01 & 500 & 2.7e-07 & 1.1e-05 \\
 & plsqr & 5.2e+02 & 2.9e+01 & 5 & 7.8e-14 & 4.9e-13 \\
\midrule
\multirow{2}{*}{\(2^{19}\)} & lsqr & - & 1.7e+02 & 500 & 9.6e-09 & 3.0e-07 \\
 & plsqr & 1.2e+03 & 6.2e+01 & 5 & 3.2e-13 & 1.1e-12 \\
\midrule
\multirow{2}{*}{\(2^{20}\)} & lsqr & - & 5.2e+02 & 500 & 9.7e-07 & 6.2e-05 \\
 & plsqr & 2.6e+03 & 1.7e+02 & 5 & 1.2e-13 & 7.6e-13 \\
\bottomrule
\end{tabular}
\caption{Time and error results of the iterative solver under the random-perturbation case, \(\alpha = 0.1\).}
\label{tab:iterative_rp_1e-1}
\end{table}

\begin{table}[tbhp]
\centering
\begin{tabular}{c cccc ccc}
\toprule
\(N\) & Method & \(t_{\pre}\) & \(t_{\iter}\) & \(n_{\iter}\) & \(r_{\s}\) & \(e_{\s}\) \\
\midrule
\multirow{2}{*}{\(2^{16}\)} & lsqr & - & 2.2e+01 & 500 & 3.0e-07 & 7.1e-06 \\
 & plsqr & 1.0e+02 & 1.5e+01 & 8 & 5.4e-13 & 2.7e-12 \\
\midrule
\multirow{2}{*}{\(2^{17}\)} & lsqr & - & 4.3e+01 & 500 & 5.1e-10 & 9.1e-09 \\
 & plsqr & 2.3e+02 & 2.6e+01 & 9 & 8.2e-15 & 3.2e-14 \\
\midrule
\multirow{2}{*}{\(2^{18}\)} & lsqr & - & 8.5e+01 & 500 & 1.2e-07 & 6.3e-06 \\
 & plsqr & 5.1e+02 & 5.2e+01 & 9 & 9.0e-14 & 5.6e-13 \\
\midrule
\multirow{2}{*}{\(2^{19}\)} & lsqr & - & 2.2e+02 & 500 & 2.1e-06 & 6.4e-04 \\
 & plsqr & 1.1e+03 & 9.4e+01 & 8 & 2.4e-13 & 2.6e-11 \\
\midrule
\multirow{2}{*}{\(2^{20}\)} & lsqr & - & 4.8e+02 & 500 & 4.3e-06 & 5.5e-04 \\
 & plsqr & 2.6e+03 & 2.4e+02 & 9 & 7.1e-14 & 1.9e-12 \\
\bottomrule
\end{tabular}
\caption{Time and error results of the iterative solver under the random-perturbation case, \(\alpha = 0.4\).}
\label{tab:iterative_rp_4e-1}
\end{table}

Figure~\ref{fig:rp_direct} and Tables~\ref{tab:direct_rp_1e-1} and~\ref{tab:direct_rp_4e-1} show the results of the direct solver under the random-perturbation case for \(\alpha = 0.1\) and \(0.4\).
First, it can be found that the relative residual and relative error are at the order of \(10^{-6}\) for all examples, which is slightly larger than the perturbation-perturbation case.
This is because the random-perturbation case is more ill-conditioned than the perturbation-perturbation case, which can also be inferred from the results of the iterative solver.
Furthermore, the runtime of the direct solver is quite similar for both the random-perturbation and perturbation-perturbation cases, which is consistent with our expectation since our algorithm is not sensitive to the distribution of sample points.

Tables~\ref{tab:iterative_rp_1e-1} and~\ref{tab:iterative_rp_4e-1} summarize the results of the iterative solver.
It is clear that lsqr fails to converge within \(500\) iterations for all examples.
On the other hand, plsqr converges within \(5\) steps for \(\alpha = 0.1\) and \(9\) steps for \(\alpha = 0.4\), which is significantly smaller than the number of iterations without preconditioning.
This demonstrates the efficiency of our method as a preconditioner for iterative methods and also implies that increasing \(\alpha\) leads to a more ill-conditioned problem, which is consistent with our intuition.
Considering the time cost, once the preconditioner is constructed~(which is the most time-consuming part), the runtime of plsqr is quite small, which is about \(170\)s for \(\alpha = 0.1\) and \(240\)s for \(\alpha = 0.4\) when \(N = 2^{20}\).

\section{Conclusions and future work} \label{sec:conclusions_future_work}

In this paper, we introduce a superfast direct solver for the least-squares problem of the type-III INUDFT problem.
Our method combines the efficient type-II INUDFT direct solver~\cite{Wilber_Epperly_Barnett_2025} and fast algorithms of HSS matrices~\cite{Levitt_Martinsson_2024, Xi_Xia_Cauley_Balakrishnan_2014}.
By approximating the type-III NUDFT matrix by a product of a NUDFT-II matrix and an HSS matrix, we develop a fast inversion algorithm, which can serve as a direct solver or an efficient preconditioner for iterative methods.
We give some theoretical results under suitable assumptions of the sample points and frequencies.
Numerical experiments are carried out to explore the properties of our fast structure and to demonstrate the efficiency of the proposed method.

Some important directions for future research include the following:
\begin{itemize}
    \item Theoretical analysis of the type-III NUDFT problem such as the properties of the NUDFT matrix.
    This would enable a much deeper understanding of the approximation.
    Current works mostly focus on the type-II case and lack a discussion of the type-III case.
    \item Investigate a method on the regularized least-squares problem.
    In many circumstances (such as when sample points and frequencies are both uniform random variables), the NUDFT matrix is ill-conditioned or even rank-deficient.
    Regularization techniques need to be added into the least-squares problem, and in this case, our method cannot be applied directly.
    The question of whether the regularized problem can be efficiently solved remains an interesting topic.
    \item Extension to high-dimensional cases.
    The challenges of the type-III INUDFT problem in high dimensions are similar to those of the type-II INUDFT.
    For two dimensions and higher dimensions, the displacement structure~\eqref{eq:nudft2_Sylvester_tildeA} of \(\Mat{\tilde{A}}\) does not hold anymore.
    In this case one can still derive some low-rank properties by the kernel function expression~\eqref{eq:nudft2_kernel_matrix}.
    However, it cannot be compressed into an HSS matrix since many blocks are not numerically low-rank.
    Designing efficient methods to solve high-dimensional INUDFT problems remains open.
\end{itemize}

\appendix

\section{The HSS rank in the type-III NUDFT approximation}

\begin{figure}[!htbp]
    \centering
    \includegraphics[width=0.5\textwidth]{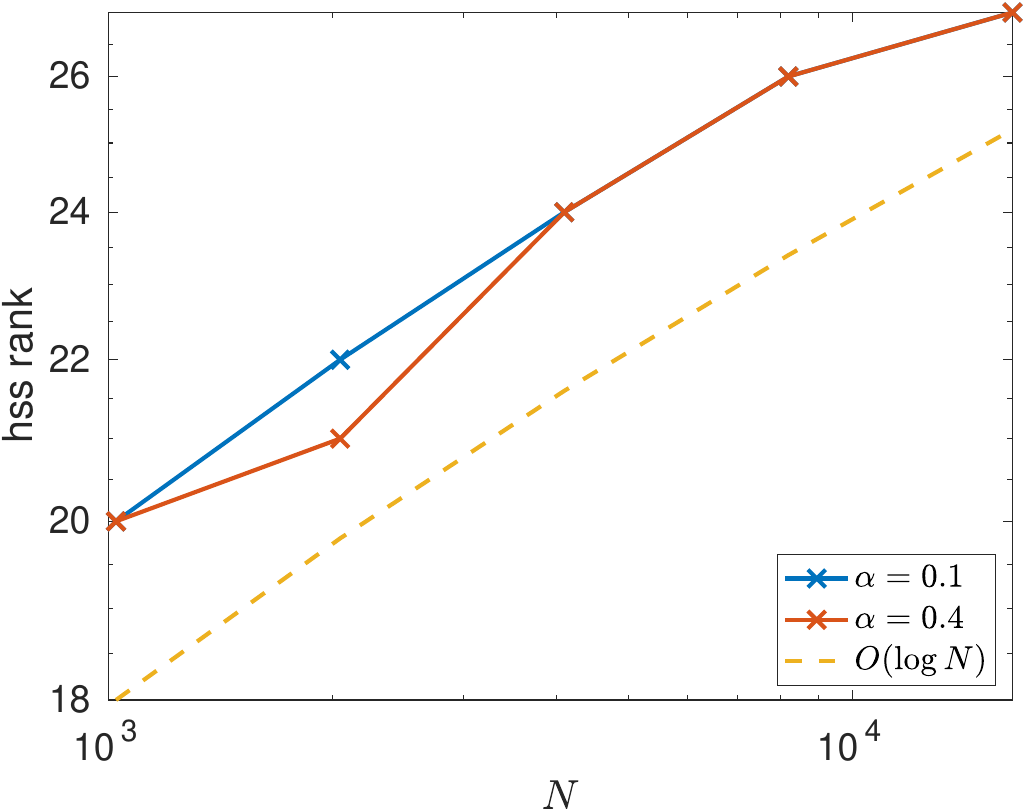}
    \caption{The \(\varepsilon\)-rank of the top-right submatrix in \(\Mat{B}_{\fast}^{\pinv} \Mat{A}\) for \(\varepsilon = 10^{-6}\).
    The size of the submatrix is \(N / 2\), where \(N = 2^p\) for \(p = 10 : 14\).}
    \label{fig:typeIII_hss_rank}
\end{figure}
We discuss the HSS rank of \(\Mat{B}_{\fast}^{\pinv} \Mat{A}\).
To do this, we extract the \(N / 2 \times N / 2\) top-right submatrix of it under the random-perturbation case and compute its \(\varepsilon\)-rank \(r_{\varepsilon}\).
The results in terms of \(\varepsilon = 10^{-6}\) are plotted in Figure~\ref{fig:typeIII_hss_rank}.
It is obvious that the numerical rank is significantly small compared to the size of the submatrix, implying the low-rank structure of \(\Mat{B}_{\fast}^{\pinv} \Mat{A}\).
Furthermore, given a fixed accuracy \(\varepsilon\), the rank of the submatrix behaves approximately as \(\bigO(\log N)\) relative to the matrix size, consistent with our intuition and discussion in section~\ref{subsec:inudft3_algorithm}.

\section{Discussion of the type-III NUDFT approximation}

In this section, we give a detailed discussion of the approximation in the fast structure~\eqref{eq:nudft3_fast_forward}, which can be viewed as an extension of section~\ref{subsec:properties_nudft3}.
Specifically, we consider the following three errors.
\begin{itemize}
    \item The relative error of the approximation:
    \begin{equation*}
        \|\Mat{A} - \Mat{B} \Mat{Q}\|_{\fro} / \|\Mat{A}\|_{\fro} \quad \text{and} \quad \|\Mat{A} - \Mat{B} \Mat{H}\|_{\fro} / \|\Mat{A}\|_{\fro}.
    \end{equation*}
    \item The relative error of the inversion:
    \begin{equation*}
        \|\Mat{A}^{\pinv} - \Mat{Q}^{\pinv} \Mat{B}^{\pinv}\|_{\fro} / \|\Mat{A}^{\pinv}\|_{\fro} \quad \text{and} \quad \|\Mat{A}^{\pinv} - \Mat{H}^{\pinv} \Mat{B}^{\pinv}\|_{\fro} / \|\Mat{A}^{\pinv}\|_{\fro}.
    \end{equation*}
    \item The relative error of the inversion on \(\range(\Mat{A})\):
    \begin{equation*}
        \|\Mat{I} - \Mat{Q}^{\pinv} \Mat{B}^{\pinv} \Mat{A}\|_{\fro} / \|\Mat{I}\|_{\fro} \quad \text{and} \quad \|\Mat{I} - \Mat{H}^{\pinv} \Mat{B}^{\pinv} \Mat{A}\|_{\fro} / \|\Mat{I}\|_{\fro}.
    \end{equation*}
\end{itemize}

\begin{table}[tbhp]
\centering
\begin{tabular}{c ccc}
\toprule
\(\sigma\) & \(1\) & \(2\) & \(3\) \\
\midrule
\(\|\Mat{A} - \Mat{B} \Mat{Q}\|_{\fro} / \|\Mat{A}\|_{\fro}\) & 7.0e-03 & 2.1e-03 & 1.7e-03 \\
\(\|\Mat{A} - \Mat{B} \Mat{H}\|_{\fro} / \|\Mat{A}\|_{\fro}\) & 5.5e-03 & 9.9e-04 & 4.9e-04 \\
\midrule
\(\|\Mat{A}^{\pinv} - \Mat{Q}^{\pinv} \Mat{B}^{\pinv}\|_{\fro} / \|\Mat{A}^{\pinv}\|_{\fro}\) & 1.4e-02 & 7.7e+01 & 5.4e+08 \\
\(\|\Mat{A}^{\pinv} - \Mat{H}^{\pinv} \Mat{B}^{\pinv}\|_{\fro} / \|\Mat{A}^{\pinv}\|_{\fro}\) & 1.1e-02 & 7.3e+01 & 1.6e+07 \\
\midrule
\(\|\Mat{I} - \Mat{Q}^{\pinv} \Mat{B}^{\pinv} \Mat{A}\|_{\fro} / \|\Mat{I}\|_{\fro}\) & 8.8e-03 & 3.8e-02 & 1.1e+05 \\
\(\|\Mat{I} - \Mat{H}^{\pinv} \Mat{B}^{\pinv} \Mat{A}\|_{\fro} / \|\Mat{I}\|_{\fro}\) & 5.9e-15 & 2.2e-13 & 8.4e-02 \\
\bottomrule
\end{tabular}
\caption{Errors of the approximations in the type-III NUDFT, \(\alpha = 0.1\).}
\label{tab:typeIII_e_sigma_1e-1}
\end{table}

\begin{table}[tbhp]
\centering
\begin{tabular}{c ccc}
\toprule
\(\sigma\) & \(1\) & \(2\) & \(3\) \\
\midrule
\(\|\Mat{A} - \Mat{B} \Mat{Q}\|_{\fro} / \|\Mat{A}\|_{\fro}\) & 2.7e-02 & 8.6e-03 & 6.6e-03 \\
\(\|\Mat{A} - \Mat{B} \Mat{H}\|_{\fro} / \|\Mat{A}\|_{\fro}\) & 2.3e-02 & 5.3e-03 & 2.4e-03 \\
\midrule
\(\|\Mat{A}^{\pinv} - \Mat{Q}^{\pinv} \Mat{B}^{\pinv}\|_{\fro} / \|\Mat{A}^{\pinv}\|_{\fro}\) & 5.9e-02 & 3.2e+01 & 2.9e+06 \\
\(\|\Mat{A}^{\pinv} - \Mat{H}^{\pinv} \Mat{B}^{\pinv}\|_{\fro} / \|\Mat{A}^{\pinv}\|_{\fro}\) & 5.4e-02 & 3.4e+01 & 4.4e+05 \\
\midrule
\(\|\Mat{I} - \Mat{Q}^{\pinv} \Mat{B}^{\pinv} \Mat{A}\|_{\fro} / \|\Mat{I}\|_{\fro}\) & 3.7e-02 & 1.4e-01 & 9.8e+03 \\
\(\|\Mat{I} - \Mat{H}^{\pinv} \Mat{B}^{\pinv} \Mat{A}\|_{\fro} / \|\Mat{I}\|_{\fro}\) & 7.2e-15 & 5.6e-13 & 3.5e-04 \\
\bottomrule
\end{tabular}
\caption{Errors of the approximations in the type-III NUDFT, \(\alpha = 0.4\).}
\label{tab:typeIII_e_sigma_4e-1}
\end{table}
Tables~\ref{tab:typeIII_e_sigma_1e-1} and~\ref{tab:typeIII_e_sigma_4e-1} summarize the results of the errors for different \(\sigma\) and \(\alpha\).
It can be found that for \(\|\Mat{A} - \Mat{B} \Mat{Q}\|_{\fro} / \|\Mat{A}\|_{\fro}\) and \(\|\Mat{A} - \Mat{B} \Mat{H}\|_{\fro} / \|\Mat{A}\|_{\fro}\), the errors decrease as \(\sigma\) increases, which is consistent with our theoretical analysis.
However, for \(\|\Mat{A}^{\pinv} - \Mat{Q}^{\pinv} \Mat{B}^{\pinv}\|_{\fro} / \|\Mat{A}^{\pinv}\|_{\fro}\) and \(\|\Mat{A}^{\pinv} - \Mat{H}^{\pinv} \Mat{B}^{\pinv}\|_{\fro} / \|\Mat{A}^{\pinv}\|_{\fro}\), the errors are less than \(1\) only when \(\sigma = 1\), meaning that we cannot use \(\sigma > 1\) for these approximations.
An increasing trend with \(\sigma\) is also observed for \(\|\Mat{I} - \Mat{Q}^{\pinv} \Mat{B}^{\pinv} \Mat{A}\|_{\fro} / \|\Mat{I}\|_{\fro}\) and \(\|\Mat{I} - \Mat{H}^{\pinv} \Mat{B}^{\pinv} \Mat{A}\|_{\fro} / \|\Mat{I}\|_{\fro}\).
Therefore we focus on the case of \(\sigma = 1\).
In this case, although \(\|\Mat{A}^{\pinv} - \Mat{Q}^{\pinv} \Mat{B}^{\pinv}\|_{\fro} / \|\Mat{A}^{\pinv}\|_{\fro}\), \(\|\Mat{A}^{\pinv} - \Mat{H}^{\pinv} \Mat{B}^{\pinv}\|_{\fro} / \|\Mat{A}^{\pinv}\|_{\fro}\), and \(\|\Mat{I} - \Mat{Q}^{\pinv} \Mat{B}^{\pinv} \Mat{A}\|_{\fro} / \|\Mat{I}\|_{\fro}\) are not very small~(at the order of \(10^{-2}\)), \(\|\Mat{I} - \Mat{H}^{\pinv} \Mat{B}^{\pinv} \Mat{A}\|_{\fro} / \|\Mat{I}\|_{\fro}\) is quite small~(at the order of \(10^{-15}\)).
It implies that when the target values are in or nearly in the range of \(\Mat{A}\), we can expect a very accurate solution by using \(\Mat{H}^{\pinv} \Mat{B}^{\pinv}\).

\section*{Acknowledgment}
The authors would like to thank the reviewers for their valuable comments and suggestions, which have greatly improved the quality of this paper.

\bibliographystyle{siamplain}
\bibliography{references}
\end{document}